%% file: NS_incompatibile.tex
\documentclass[12pt,a4paper]{article}

\input{ex_shared}

\begin{document}

\maketitle

\begin{abstract}
  In this paper we shall consider the Navier-Stokes equations in the half plane with Euler-type initial conditions, 
  i.e. initial conditions which have a non-zero tangential component at the boundary. 
  Under analyticity assumptions for the data, we shall prove that the solution exists for a small time independent of the viscosity. 
  The solution is constructed through a composite asymptotic expansion involving the solutions of the Euler and Prandtl equations, plus an error term.
  The norm of the error goes to zero with the square root of the viscosity. 
  The Prandtl solution contains a singular term, which influences the regularity of the error.  
  The error term is written as the sum of a first order Euler correction, a first order Prandtl correction, and a further error term. 
  The use of an analytic setting is mainly due to the Prandtl equation.
\end{abstract}



\section{Introduction}
In this paper, we shall study the solutions of the 3D Navier-Stokes(NS)  equations in the half-space 
when the initial datum and the boundary datum are incompatible. 
At the boundary, we shall impose  the no-slip boundary condition and, therefore, 
the incompatibility means that we shall consider initial data having non-zero tangential component at the boundary;
 we shall call these data {\it Euler type} or {\it not well prepared} initial data. 
We shall establish the existence and uniqueness of the NS solutions for a  time short but independent of the viscosity. 
The central hypothesis we shall impose is the analyticity of the initial datum. 

The study of the NS equations in the half-space (or half-plane in 2D) is a central subject in the mathematical theory of fluid dynamics 
because it is a prototypical case where one can study the interaction between a fluid and a wall; 
and explaining such interactions and the resulting phenomena (like boundary layer formation, vorticity generation, transitions) 
motivates a significant part of the interest in fluids. 
The introduction by S. Ukai of an exact formula for solutions to the Stokes equations in the half-space \cite{Ukai1987}, 
had a significant impact on the mathematical theory of the NS equations; 
since Ukai's result, several papers appeared concerning the existence and uniqueness of solutions for the Navier-Stokes 
equations in the half-space with general initial data, see for example \cite{CPS2000, Mar2008}  and references therein. 
However, these results, by one side, do not allow Euler type initial data, i.e., data with a non zero tangential component at the boundary; 
on the other hand, they rely on estimates that, in the zero viscosity limit, would degrade. 

Euler type initial data are important classical configurations; among them, we mention the impulsively started disk or plate \cite{Batc} 
and the flow generated by the interaction between a wall and a point vortex or a core vortex \cite{Lamb}. 
Incompatible data are an interesting subject also from a numerical point of view; in \cite{BF1999,  CQT2010,CQT2011,HTZ2013,Tem2006}, 
one can find in-depth studies of how the lack of compatibility between initial and boundary data can lead to loss of numerical accuracy;
and on the appropriate compatibility conditions, one should impose to ensure the required accuracy. 
In incompressible fluid dynamics Euler type data typically arise when, to the Navier-Stokes equations, one imposes  
initial data that are stationary solutions of the Euler equations and are therefore of great interest for the numerical and 
theoretical study of the boundary layer theory; 
for example, they have been used to test the Prandtl equations' effectiveness to reproduce the separation phenomena 
occurring at the boundary, see \cite{GSS2011,GSSC2014,OC2002}.

The study of the zero viscosity limit of the incompressible Navier-Stokes solutions is one of the most important and challenging subjects 
in the mathematical theory of fluid dynamics. 
Formally, one should have the convergence of the NS solutions to the Euler solutions. 
However, NS equations are a singular perturbation of the Euler equations and, from the physical point of view, 
the interaction of the flow with the boundary creates sharp gradients in the normal directions;  
this situation requires the introduction of a boundary layer corrector to correct the discrepancy between the 
no-slip boundary condition and the Euler solution: classically, such a corrector is the Prandtl equations solution. 
The lack of a completely satisfactory theory of the Prandtl boundary layer equations has been, so far, 
a major obstacle in the characterization of the zero viscosity solutions of the NS equations. 
The monotonic setting \cite{OS1997} (through the use of the Crocco transformation as in the original Oleinik's work, or through the 
use of energy estimates, as in \cite{AWXY2015,LWY2017,MW2015}),  
and the analytic setting \cite{Asa1988,CLS2013,KMVW2014,KV2013,LCS2003,SC1998A}, 
are two instances where one can control the ill-posedness and instability phenomena recently discovered 
(\cite{Gre2000,G-VD2010,GN2011, G-VN2012,LWY2016,LY2017}); 
so that it is possible to prove the well-posedness of Prandtl's equations. 
More recently, similar results have been achieved in the Gevrey class (rather than analytic) data 
\cite{G-VM2015,LY2020} assuming some structural hypothesis on the initial data (like the non degeneracy of 
the singular point) or without such assumptions in \cite{DG-V2019}. 

The construction, see \cite{SC1998A}, of the analytic Prandtl solution allowed to show the convergence of the NS solution 
in the half-plane as well in the half-space.  
In \cite{SC1998B} the authors constructed the NS solution through an asymptotic matching between the Prandtl and the Euler solutions, 
plus a correction which was shown to be $O(\sqrt{\nu})$. 
From which the convergence to Euler and Prandtl solutions in their respective domain of validity followed. 
See also \cite{CS1997}  for a different geometry (the outside of a disk) and, more 
recently, the papers  \cite{WWZ2017}, where the above results have been derived by the use of energy methods, and  
\cite{NN2018}, where the inviscid limit was derived without the use of a boundary-layer corrector. 

Recently, in the 2D half-plane, all the above has been improved, showing that analyticity is necessary only close to the boundary; 
this important result has been achieved first by Maekawa in \cite{Mae2014}, where zero vorticity close to the boundary was assumed;  
this result has been generalized to the half-space in \cite{FTZ2018}. 
Kukavica, Vicol, and Wang \cite{KVW2020}, in the half-plane, 
proved that data analytic only near the boundary and Sobolev regular elsewhere are enough to obtain the validity of the inviscid limit: these results have been obtained through energy methods applied to the vorticity formulation 
of the Navier-Stokes equations, without the use of the matched asymptotic expansion of the solution. 
This result has been extended to a 3D setting in \cite{Wan2020}.

In 1994 Kato initiated a different line of research. In \cite{Kato84}, he showed that the convergence of the NS solutions to the 
Euler solution is equivalent to the vanishing of the dissipation in a sublayer of thickness $O(\nu)$, therefore smaller than
the Prandtl layer.  
Temam and Wang, in \cite{TW1997}, see also \cite{CW2007,Wang2001}, obtained a condition leading to convergence, 
based on the pressure gradient at the boundary, which is interesting because the appearance of high stress at the boundary 
is the first indicator of the deviation of the Navier-Stokes solutions from the Prandtl solution and the precursor of vortices 
formation and subsequent separation, see \cite{GSS2011,GSSC2014,OC2002}. 
During the last 15 years, the Kato criterion has been improved, interpreted in terms of different quantities of physical interest, e.g., vorticity \cite{Kel2008,BT2013,Kel2017} or tangential velocity and velocity gradient at the boundary \cite{CEIV2017}, 
combined with Oleinik's monotonicity setting \cite{CKV2015}. 
Recently, in  \cite{CV2018,CLNV2019,DN2019} the authors have derived interesting criteria based on the vorticity in the interior of the domain; 
being impossible to report all the implication of this line of research, we refer the interested reader to \cite{BT2013} or to the recent review 
paper \cite{MM2016}. 

Further cases where it has been possible to establish the inviscid limit are flows with symmetries \cite{BW2002,HMNW2012,Kel2009,LMN2008,LMNT2008}, 
flows with anisotropic viscosity  \cite{Mas1998}. 

As we mentioned before, we want to construct the Navier-Stokes equation solution starting from an initial datum that has a non-zero tangential slip. 
Moreover, we want that the time of existence of the solution does not shrink to zero when the viscosity goes to zero; 
this is not trivial because it is well known that, in 3D, the time of existence of the NS solutions is short and that, in the presence of boundaries, 
the existence time, in general, would go to zero with the viscosity.  
Therefore, we shall have to handle both the initial layer, due to the presence of the initial discontinuity, and the boundary layer, 
due to the mismatch between the no-slip boundary condition and the Euler solution. 
Following \cite{SC1998B}, we shall decompose the solution to the Navier-Stokes equations 
as the sum of the solution to the Euler equations, the solution to the Prandtl equations with non-compatible data,
and a remainder; due to the presence of the initial layer, we shall see that further decomposition of the remainder is necessary. 
In our procedure we shall rely on the well-posedness result for Prandtl's equations with non-compatible data \cite{CLS2013}. 

We also point out the recent work that has tackled similar problems, although considering a linearized version of the NS equations
\cite{Gie2013,GKM2018}.

\subsection{Asymptotic decomposition of the NS solution}

Consider the Navier-Stokes equations in the half-plane or in the half-space $\Pi^+=\{(x,y): y \geq 0\}$ where, in 3D, $x=(x_1,x_2)$:
 \begin{equation} \label{Navier-Stokes}
\begin{split}
\partial_t \textbf{u} - \nu \Delta \textbf{u} + \textbf{u} \cdot \nabla \textbf{u} + \nabla p =0,  \\
\nabla \cdot \textbf{u} =0, \\
\textbf{u} \big|_{t=0}=\textbf{u}_0, \\
\gamma \textbf{u}= 0. 
\end{split}
\end{equation}
In the above equation $\gamma$ is the trace operator, i.e., for regular functions, 
\begin{equation}
\gamma \textbf{v}= \textbf{v}(y=0).
\end{equation}

When the viscosity  $\nu$  is small, its effects are mainly concentrated in a layer near the boundary, whose thickness is 
proportional to $\varepsilon=\sqrt{\nu}$. 
For this kind of singularly perturbed problems, one looks for solutions that, at the formal level, can be written as \cite{VanDyke64}:
\begin{equation}
\textbf{u}=\textbf{u}^{out}_{(0)}+\textbf{u}^{inn}_{(0)}+\varepsilon\left(\textbf{u}^{out}_{(1)}+\textbf{u}^{inn}_{(1)}\right) + 
....+ \varepsilon^N \left(\textbf{u}^{out}_{(N)}+\textbf{u}^{in}_{(N)}\right) +O(\varepsilon^{N+1})  ,   \label{expansionVD}
\end{equation}
having carried out the expansion up to order $N$. 
In the above expansion  $\textbf{u}^{out}_{(i)}(x,y,t)$ and $\textbf{u}^{in}_{(i)}(x,Y,t)$ for  $i=0,...N$ describe the outer flow (away from the boundary)
and the inner flow (close to the boundary) respectively, being $Y=y/\varepsilon$ the rescaled variable. 

The rigorous proof that an expansion like \eqref{expansionVD} holds in the compatible case, was given in \cite{SC1998B}, where the authors proved that 
the solution of the NS equations can be written as: 
$$
\textbf{u}= \textbf{u}^E+ \bar{\textbf{u}}^P+\varepsilon\left( \textbf{u}^E_1+ \bar{\textbf{u}}^P_1 \right)+\varepsilon  \textbf{e}. 
$$
In the above expansion  the outer leading order term, $\textbf{u}^E $, is the solution of the Euler equations: 
 \begin{equation} \label{Euler}
\begin{split}
\partial_t \textbf{u}^E + \textbf{u}^E \cdot \nabla \textbf{u}^E + \nabla p^E =0,  \\
\nabla \cdot \textbf{u}^E =0, \\
\textbf{u}^E \big|_{t=0}=\textbf{u}_0, \\
\gamma_n\textbf{u}^E=0, \\
\end{split}
\end{equation}
where $\gamma_n\textbf{v}$ is the normal component of the trace of $\textbf{v}$. 

The  leading order inner solution  $\bar{\textbf{u}}^P$, is linked to the solution of the Prandtl equations 
as follows: let $u^P$ the solution of Prandtl's equations 
\begin{equation} \label{Prandtl}
\begin{split}
(\partial_t-\partial_{YY})u^P + u^P \partial_x u^P+v^P\partial_Yu^P- (\partial_t+u^E\big|_{y=0} \partial_x)u^E\big|_{y=0}=0, \\
u^P\Bigr|_{Y=0}=0, \, \, \, \, \, u^P\Bigr|_{Y\rightarrow +\infty}=u^E(x,y=0,t), \\ 
u^P\Bigr|_{t=0}=u_0(x,y=0),
\end{split}
\end{equation}
where \begin{equation}
v^P= - \int \limits_0^Y \partial_xu^P d Y'. 
\end{equation}
We define the inner solution as:
\begin{equation}
\bar{\textbf{u}}^P=(\tilde{u}^P, \varepsilon \bar{v}^P) ,  \label{uboldbar}
\end{equation}
where
\begin{equation}\label{utildeP}
\tilde{u}^P=u^P-u^E\big|_{y=0},
\end{equation}
\begin{equation}\label{vbarP}
\bar{v}^P= \int \limits_Y^{+\infty} \partial_x \tilde{u}^P d Y'.
\end{equation}
The $\bar{\textbf{u}}^P$ defined above decays away from the boundary, so that it does not interfere with the outer solution $\textbf{u}^E$. 
Moreover $\bar{\textbf{u}}^P$ cancels the tangential flow generated by  $\textbf{u}^E$. 
However $\bar{\textbf{u}}^P$ generates a normal inflow at the boundary $Y=0$, see \eqref{vbarP}. 
This inflow is canceled by the first order outer solution $\varepsilon\textbf{u}^E_1$. 
This outer solution generates tangential flow which is canceled by  first order inner solution $\varepsilon\bar{\textbf{u}}^P_1$.

The overall remainder $\varepsilon  \textbf{e}$ closes the procedure. 
Notice, however, that the remainder, which at the formal level should be $O(\varepsilon^2)$, 
due to the nonlinear interactions between the boundary layer terms and the outer Euler terms, 
can be proven to be $O(\varepsilon)$ only. 

\subsection{Statement of the main result} \label{informal_main_result}
We give an informal statement of the main result of the present paper.
\begin{theorem}[Informal statement]
Suppose that $\textbf{u}_0$ has analytic regularity, with $\gamma_n \textbf{u}_0=0$. Then an analytic solution to the Navier-Stokes equations exists for a time independent of the viscosity; this solution can be written as
\begin{equation}\label{informalmain}
\textbf{u}=\textbf{u}^{out}_{(0)}+\textbf{u}^{inn}_{(0)}+\varepsilon\left(\textbf{u}^{out}_{(1)}+\textbf{u}^{inn}_{(1)}+ \textbf{e}\right) ,
\end{equation}
where the outer solutions $\textbf{u}^{out}_{(0)}$ and $\textbf{u}^{out}_{(1)}$ depend on $(x,y,t)$, while the inner solutions 
$\textbf{u}^{inn}_{(0)}$ and $\textbf{u}^{inn}_{(1)}$ depend on $(x,Y,t)$, and are exponentially decaying outside the boundary layer. 

Moreover $\textbf{u}^{inn}_{(0)}$, $\textbf{u}^{out}_{(1)}$ and $\textbf{u}^{inn}_{(1)}$ have singular time derivatives at the initial time. 
\end{theorem}

The rest of this subsection is devoted to illustrating the physical meaning of the terms appearing in the asymptotic expansion of the NS solution, 
equation \eqref{informalmain}. 

\subsubsection{The leading-order outer solution}
The term $\textbf{u}^{out}_{(0)}(x,y,t)$ is the solutions $\textbf{u}^E$ of the Euler equations \cref{Euler}. 
The existence of an analytic solution for the Euler equation has been established in the literature, see among others \cite{BB77,LO1997,SC1998A}; 
therefore, in constructing this part of the NS solution, we shall refer to these results. 

\subsubsection{The leading-order inner solution}

The term  $\textbf{u}^{inn}_{(0)}=\bar{\textbf{u}}^P(x,Y,t)$ correlates to the solution of the Prandtl equations \eqref{Prandtl}
through subtraction of a constant (in $Y$) function that makes it decaying outside the BL. 
The incompatibility,  in equations \eqref{Prandtl}, between the initial data and the boundary data leads to the occurrence of a singular part  $\textbf{u}^S$, 
so that the  $\bar{\textbf{u}}^P$ must be  decomposed in a regular and a singular part: 
 \begin{equation}
 \bar{\textbf{u}}^P= \bar{\textbf{u}}^R+\textbf{u}^S. \label{prandtl_decomposition}
 \end{equation}
The singular part  $\textbf{u}^S=(u^S,\varepsilon v^S)$ is in the form of an initial layer;  $u^S$ solves the following system:
 \begin{equation}
\begin{split}
(\partial_t - \partial_{YY})u^S=0, \\
u^S\bigr|_{Y=0}= - u_{00}(x), \\
u^S|_{t=0}=0,
\end{split}
\end{equation}
where $u_{00}$ is the trace at the boundary of the tangential part of the initial NS datum, i.e., 
\begin{equation}
u_{00}=\gamma u_0; 
\end{equation}
clearly, were the initial and boundary data compatible, $u_{00}=0$ and, consequently, $u^S=0$.
One can write the explicit expression of $u^S$, 
\begin{equation}
u^S=-\gamma u_0 \int \limits_Y^{+\infty} \frac{e^{-\frac{\sigma^2}{4t}}}{\sqrt{\pi t}} d \sigma, \label{u_singular}
\end{equation}
while $v^S$ derives from the incompressibility condition, 
\begin{equation}
v^S= \int \limits_Y^{+\infty} \partial_x u^S d Y'.\label{v_singular}
\end{equation}
Notice how the tangential part $u^S$ has $Y$- and $t$- derivatives that, at the initial time, are singular at the boundary. 
For the rest of the paper, we shall say that a function has a "gaussian" singularity if it is bounded away from the origin $(Y,t)=(0,0)$ and behaves near the origin like \begin{equation}
t^{-\beta} p(Y/\sqrt{t}) e^{-Y^2/{c t}}
\end{equation} for some $\beta >0$, $c>0$ and some polynomial $p$. This kind of singularity is concentrated at the boundary at the initial time, and that can be tamed by multiplying the function by suitable powers of $Y$ and/or $t$.
The regular part $\bar{\textbf{u}}^R$ 
can be proven to be more regular. 
The fact that the Prandtl equation with incompatible data admits the solution given in \eqref{prandtl_decomposition} was proven in  \cite{CLS2013}. 

In the sequel it will be useful to introduce the following notation for the value at the boundary of the influx generated by the Prandtl solution, 
 \begin{equation}\label{VPrandtlboundary}
 	\begin{split}
	g\equiv  \gamma \bar{v}^P= - \int \limits_0^{+ \infty}  \partial_x \tilde{u}^P d Y'= g^R+g^S, \\
	 g^R=- \int \limits_0^{+ \infty}  \partial_x \tilde{u}^R d Y'  \qquad g^S=- \int \limits_0^{+ \infty}  \partial_x u^S d Y'. =-\sqrt{\frac{4}{\pi}}\sqrt{t}\partial_x u_0(x,y=0)
\end{split}
\end{equation}

\subsubsection{The first-order correction to the outer solution}

The first order outer flow $\textbf{u}^{out}_{(1)}=\textbf{u}^E_{(1)}(x,y,t)$ is the solution of a linearized Euler's system: 
\begin{equation}\label{CorrEul}
	\begin{split}
		\partial_t \textbf{u}^E_{(1)} + \textbf{u}^E_{(1)} \cdot \nabla \textbf{u}^E + \textbf{u}^E \cdot \nabla \textbf{u}^E_{(1)} + \nabla p^{E}_{(1)} = \textbf{0}, \\
		\nabla \cdot \textbf{u}^E_{(1)}=0, \\
		\gamma_n \textbf{u}^E_{(1)}= -\gamma \bar{v}^P ,\\
		\textbf{u}^E_{(1)}\big|_{t=0}= \textbf{0}. 
	\end{split}
\end{equation}
The role of the correction to the Euler flow is to cancel the normal inflow generated at the boundary by the boundary layer corrector $\bar{\textbf{u}}^P$. 
We shall decompose $\textbf{u}^E_{(1)}$ in a regular and a singular part:
\begin{equation}
\textbf{u}^E_{(1)}=\boldsymbol{w}^R+\boldsymbol{w}^S. \label{euler_first_decomposition}
\end{equation}
We shall see that the singular part can be written as $\sqrt{t}\boldsymbol{w}^S_b(x,y)$, where $\boldsymbol{w}^S_b$ does not depend on $t$. This time behavior means that the singularity, formed at $t=0$ in the inner flow $\textbf{u}^{inn}_{(0)}$,  does not remain confined in the boundary layer, but instantaneously propagates  in the whole space with an $O(\varepsilon)$ intensity. 
This is consistent with the parabolic nature of the Navier-Stokes equations, leading to an infinite speed of propagating disturbances. 

While the time derivative of $\textbf{u}^E_{(1)}$ is singular everywhere at the initial time, its growth as $t$ goes to zero is less pronounced than the one of $\partial_t u^S$: the reason is that the singularity has been passed at $O(\varepsilon)$ terms through the incompressibility condition that gives $g^S$ in \eqref{VPrandtlboundary}, and the integration in $Y$ has a regularizing effect on "gaussian" singularities. For the same reason, we shall see that the overall error $\textbf{e}$ is regular, with a bounded time derivative.

\subsubsection{The first-order correction to the inner solution}

The first order inner flow  $\textbf{u}^{inn}_{(1)}=\bar{\textbf{u}}^P_{(1)}(x,Y,t)=\left(\bar{u}^P_{(1)}, \varepsilon \bar{v}^P_{(1)}\right) $ is the solution of a 
linearized Prandtl's system, 
\begin{equation}\label{uP1syst}
	\begin{split}
		(\partial_t - \partial_{YY}) \bar{u}^P_{(1)}=0, \\
		\gamma \bar{u}^P_{(1)}=-\gamma u^E_{(1)}, \\
		 \bar{u}^P_{(1)}(t=0)=0, %\\
%		\bar{v}^P_{(1)}=  \int \limits_Y^{+\infty} \partial_x \bar{u}^P_{(1)} d Y',
	\end{split}
\end{equation}
with the normal component given by the incompressibility condition,
\begin{equation} \label{vP1}
		\bar{v}^P_{(1)}=  \int \limits_Y^{+\infty} \partial_x \bar{u}^P_{(1)} d Y'.
\end{equation}
The role of $\bar{\textbf{u}}^P_{(1)}$ is to correct the tangential component at the boundary generated by ${\textbf{u}}^E_{(1)}$. 
The term $\bar{\textbf{u}}^P_{(1)}$ admits a decomposition in a regular part and a singular part with a $gaussian$ singularity
less severe then the one in $\textbf{u}^{inn}_{(0)}$. 

\subsubsection{The error $\textbf{e}(x,Y,t)$}

The term $\textbf{e}(x,Y,t)$ is an overall error that closes the asymptotic expansion. 
It satisfies a NS type equation with a source term and with a boundary condition that cancels the inflow generated by $\bar{\textbf{u}}^P_{(1)}$, 
without generating tangential flow:
\begin{equation}\label{corrector_e}
	\begin{split}
		\partial_t \textbf{e}+ \left(\textbf{u}^{NS}_{(0)}+\varepsilon \textbf{u}^{NS}_{(1)} \right) \cdot \nabla \textbf{e} + 
\textbf{e}\cdot \nabla  \left(\textbf{u}^{NS}_{(0)}+\varepsilon \textbf{u}^{NS}_{(1)} \right)  +\varepsilon \textbf{e}\cdot \nabla \textbf{e}+
\nabla p^e =\boldsymbol{\Xi}, \\
\nabla\cdot \textbf{e}=0, \\
\gamma  \textbf{e}=\left(0,-\gamma \bar{v}^P_{(1)}\right), \\
 \textbf{e}(t=0)=0, 
		\end{split}
	\end{equation}
where we have defined the zero-th and first order approximation of the NS solution:
$$
 \textbf{u}^{NS}_{(0)}=\textbf{u}^E+\bar{\textbf{u}}^P, \qquad \textbf{u}^{NS}_{(1)}=\textbf{u}^E_{(1)}+\bar{\textbf{u}}^P_{(1)}. 
$$
The source term $\boldsymbol{\Xi}$ is generated by the discrepancies between the NS equation and the equations 
satisfied by the approximating terms in the asymptotic expansion. The explicit expression can be found in the Appendix.

A key point in the proof of the regularity of $\textbf{e}$ is that, in its equations, the singular terms multiply terms which go to zero with $t$ or $Y$ in a sufficiently fast way. 

All the singular terms appearing in \eqref{informalmain} have explicit expressions and finite $L^2$ norms in appropriate spaces of holomorphic functions. 
The expansion \eqref{informalmain} implies in particular the validity of the inviscid limit in the energy norm with an $O(\sqrt{\varepsilon})$ rate, since $||\bar{\textbf{u}}^P||_{L^2_{xy}}= \sqrt{\varepsilon}||\bar{\textbf{u}}^P||_{L^2_{xY}}$.

\subsubsection{Plan of the paper}

The organization of the paper is the following: in  \cref{sezionespazifunzionali} we introduce the function spaces needed to prove the validity of the inviscid limit, and we analyze some of their properties.  In  \cref{sezioneteoremadicauchykow} we introduce the abstract Cauchy-Kowalewski theorem; this theorem is essentially a fixed point method used to prove the existence of solutions to differential equations in a scale of Banach spaces. In  \cref{sezionerisultatoprincipale} we present the main result of the paper. The other sections deal with the terms of order one and the remainder deriving from the asymptotic expansion of $\textbf{u}^{NS}$.

\section{Function spaces}\label{sezionespazifunzionali}
Define the strip $D(\rho)$, the angular sector $\Sigma(\theta)$ and the conoid $\Sigma(\theta,a)$ as follows \begin{equation}
D(\rho)= \{ x \in \mathbb{C}: |Im(x)|<\rho \},
\end{equation}
\begin{equation}
\Sigma(\theta)= \{ y \in \mathbb{C}: Re(y) > 0, \, |Im(y)| < Re(y) \tan( \theta) \},
\end{equation}
\begin{equation}
\begin{split}
\Sigma(\theta,a)= \{ y \in \mathbb{C}: 0< Re(y) \leq a, \, |Im(y)| < Re(y) \tan( \theta) \} \\ \cup \{y \in \mathbb{C}: Re(y)>a, \, |Im(y)| < a \tan(\theta) \}.
\end{split}
\end{equation}
In what follows, we shall deal with spaces of functions analytic in some of the above domains; 
in those spaces, the chosen paths of integration in the $x$ variable are lines parallel to the real axis, 
while we shall adopt piecewise linear paths for the $y$ and $Y$ variables.  
\begin{equation}
\Gamma(b)= \{ x \in \mathbb{C}: Im(x)=b \},
\end{equation}
\begin{equation}
\begin{split}
\Gamma(\theta',a)= \{ y \in \mathbb{C}: 0<Re(y) \leq a, \, Im(y)= Re(y) \tan(\theta') \} \cup \\ \{y \in \mathbb{C}: Re(y)>a, \, Im(y)=a \tan(\theta') \}.
\end{split}
\end{equation}
\begin{definition}
The space $H^{l,\rho}$ is the set of all complex functions $f(x)$ such that \begin{enumerate}
\item f is analytic in $D(\rho)$;
\item $|f|_{l,\rho}= \sum \limits_{|\alpha| \leq l} \sup \limits_{|\lambda|<\rho} |\partial_x^\alpha f(\cdot + i \lambda)|_{L^2(\mathbb{R})} < \infty$.
\end{enumerate}
\end{definition} 
When dealing with the 3D Navier Stokes system, $x$ is a 2D vector, and $\alpha$ is a multi-index, while in the two-dimensional case, 
$x$ is a scalar. 
%Essentially, $H^{l,\rho}$ spaces are Sobolev spaces which use, for every derivative, Hardy spaces instead of $L^2$ spaces. 
The use of an $L^2$ norm instead of a generic $L^p$ norm allows an important characterization for Hardy spaces, see \cite{PW1934} for a proof:
\begin{theorem}[Paley-Wiener Theorem for the strip] 
Let $f \in L^2(\mathbb{R})$, $\rho>0$, denote with $\hat{f}$ the Fourier transform of $f$. The following are equivalent: \begin{enumerate}
\item f is the restriction to the real line of a function holomorphic on the strip $D(\rho)$, with $\sup \limits_{|\lambda|<\rho} |f(\cdot+i\lambda)|_{L^2(\mathbb{R})}<\infty$;
\item $e^{\rho |\xi|} \hat{f} \in L^2(\mathbb{R})$.
\end{enumerate}
\end{theorem}
We shall use extensively this result, and we shall often work in the Fourier variable $\xi'$ corresponding to the physical variable $x$: to simplify the notation, we denote with the same symbol a function $f$ and its Fourier transform with respect to the $x$ variable, and similarly we use the same notation for a pseudodifferential operator and its symbol.

\subsection{Function spaces for the outer flow:  zero-th and first-order Euler equations}

In this section, we shall define the appropriate function spaces for studying Euler equations and the first order correction to the Euler equations. 
First, we introduce the space of functions, depending on $x$ and $y$. 
To construct the solution of the Euler equations,  the main tool is the half-plane Leray projector, which allows as many derivatives in $x$ and $y$.  
This reflects in the following definition. 
\begin{definition}
$H^{l,\rho,\theta}$ is the set of all functions $f(x,y)$ such that: \begin{enumerate}
\item $f$ is analytic in $D(\rho) \times \Sigma(\theta,a)$;
\item $|f|_{l,\rho,\theta}= \sum \limits_{|\alpha_1|+\alpha_2 \leq l} \sup \limits_{|\theta'|<\theta} ||\partial_y^{\alpha_2} \partial_x^{\alpha_1} f|_{0,\rho}|_{L^2(\Gamma(\theta',a))}< \infty$.
\end{enumerate}
\end{definition}
The initial value for the Euler equations is given in $H^{l,\rho,\theta}$, with $l\geq 6> 1 + d/2$ and $\theta< \pi/4$.

We now pass to introduce function spaces with time dependence. 
In all the spaces defined below, the width of the analyticity domain diminishes linearly with the time $t$. 

For a given Banach scale $ \{ X_{\rho} \}_{0<\rho \leq \rho_0}$, with $X_{\rho''} \subset X_{\rho'}$ and $|\cdot|_{\rho'} \leq |\cdot|_{\rho''}$ when $\rho' \leq \rho'' \leq \rho_0$, denote with $B^j_{\beta}([0,T],X_{\rho_0})$ the set of all functions $f$ such that, for $k=0,...,j$, $\partial_t^k f$ is continuous from $[0,\tau]$ to $X_{\rho_0-\beta \tau}$ $\forall 0< \tau \leq T \leq \rho_0/\beta$, with norm \begin{equation}
|f|_{k,\rho,\beta}= \sum \limits_{j=0}^{k} \sup \limits_{t \in [0,T]} |\partial_t^j f(t)|_{\rho_0- \beta t}.
\end{equation}

The following spaces are where one can prove the existence of the outer solutions.
\begin{definition}
$H^{l,\rho,\theta}_{\beta,T}$, with $T \leq \min \{\rho/\beta, \theta/\beta\}$ is the space of all functions $f(x,y,t)$ such that $ f \in \bigcap \limits_{i=0}^{l} B_\beta^{i}([0,T], H^{l-i,\rho,\theta})$, with norm \begin{equation}
|f|_{l,\rho,\theta,\beta,T}= \sum \limits_{i=0}^{l} \sup \limits_{t \in [0,T]} |\partial_t^i f(\cdot,\cdot,t)|_{l-i,\rho-\beta t, \theta - \beta t}.
\end{equation}
\end{definition}
The above space is the natural space for analytic solutions of the Euler equations; see Theorem 4.1 in \cite{SC1998A}. 
\begin{definition}
$\tilde{H}^{l,\rho,\theta}_{\beta,T,1}$ is the set of all functions $f(t,x,y)$ such that \begin{enumerate}
\item $f \in B^0_\beta([0,T], H^{l,\rho,\theta}) \cap B^1_\beta([0,T],H^{l-1,\rho,\theta})$ ;
\item $|f|_{l,\rho,\theta,\beta,T,1}= \sup \limits_{t \in [0,T]} |f|_{l,\rho-\beta t,\theta-\beta t} + 
\sup \limits_{t \in [0,T]}  |\partial_t f|_{l-1,\rho-\beta t,\theta-\beta t} < \infty$.
\end{enumerate}
\end{definition}
The above space is the natural space where to prove the existence of the correction of Euler flow, $\textbf{u}^E_{(1)}$:  
the singularity deriving from Prandtl equations implies that only one time-derivative is allowed; 
on the other hand, the regularity in $y$ depends only on operators (like the Leray projector) 
with some symmetry between the behavior in $x$ and $y$, 
so that we can use as many derivatives with respect to $y$ as for $x$.

\subsection{Function spaces for the inner flow:  zero-th and first order boundary layer equations}

We now introduce the  space of functions analytic in $x$ and $Y$, in the strip $D(\rho)$ and in the cone $\Sigma(\theta)$, respectively. 
Moreover, we impose exponential decay in the $Y$-variable so that for large $Y$ (away from the boundary, in the outer flow), the boundary layer-type 
solutions do not influence the Euler solutions.  

\begin{definition}
The space $K^{l,\rho,\theta,\mu}$ with $\mu >0$ is the set of all functions $f(x,Y)$ such that \begin{enumerate}
\item f is analytic in $D(\rho) \times \Sigma(\theta)$;
\item $\partial_Y^{\alpha_1} \partial_x^{\alpha_2} f \in C^0(\Sigma(\theta), H^{0,\rho})$ for any $\alpha_1 \leq 2$ and $|\alpha_2| \leq l - \alpha_1$;
\item $|f|_{l,\rho,\theta,\mu}= \sum \limits_{\alpha_1 \leq 2} \sum \limits_{|\alpha_2| \leq l - \alpha_1} \sup \limits_{Y \in \Sigma(\theta)} e^{\mu Re(Y)} | \partial_Y^{\alpha_1} \partial_x^{\alpha_2} f(\cdot,Y)|_{0,\rho} < \infty$.
\end{enumerate}
\end{definition}
The asymmetry on the regularity requirements in $x$ and $Y$ is due to the diffusion in $Y$ that 
allows regularity only up to second-order $Y$-derivatives (unless one imposes stronger compatibility conditions). 

We now introduce the functions depending on $x$ and $t$. The presence of the diffusion allows regularity in $t$ only up to first-order derivative. 
\begin{definition}
$K^{l,\rho}_{\beta,T}$ is the set of all functions $f(x,t) \in B^0_\beta([0,T], H^{l,\rho}) \cap$ \\$  B^1_\beta([0,T],H^{l-1,\rho})$, with norm \begin{equation}
|f|_{l,\rho,\beta,T}= \sum \limits_{i=0}^{1} \sum \limits_{|\alpha|\leq l} \sup \limits_{0 \leq t \leq T} |\partial_t^i \partial_x^{\alpha} f(\cdot,t)|_{0,\rho-\beta t}.
\end{equation}
\end{definition}

The next space is where we shall prove the existence of the regular part of the solutions of the zeroth and first-order BL equations. 
The diffusion in $Y$ allows regularity only up to second-order $Y$-derivatives and first-order $t$-derivatives. 
\begin{definition}
The space $K^{l,\rho,\theta,\mu}_{\beta,T}$ is the set of all functions $f(x,Y,t)$ such that: \begin{enumerate}
\item $f \in C^0([0,t], K^{l,\rho-\beta t,\theta- \beta t,\mu- \beta t})$ and $\partial_t \partial_x^{\alpha} f \in C^0([0,t], K^{0,\rho-\beta t,\theta- \beta t,\mu- \beta t})$, with $t \leq T \leq \min \{\rho/\beta,\theta/\beta,\mu/\beta \}$ and $|\alpha| \leq l-2$;
\item $|f|_{l,\rho,\theta,\mu,\beta,T} = \sum \limits_{\alpha_1 \leq 2} \sum \limits_{ \alpha_1 + |\alpha_2| \leq l} \sup \limits_{0 \leq t \leq T} |\partial_Y^{\alpha_1} \partial_x^{\alpha_2} f(\cdot,\cdot,t)|_{0,\rho- \beta t, \theta - \beta t, \mu - \beta t} + \\ \sum \limits_{|\alpha| \leq l-2} \sup \limits_{ t \in [0,T]} |\partial_t \partial_x^{\alpha} f(\cdot,\cdot,t)|_{0,\rho- \beta t, \theta - \beta t, \mu - \beta t} < \infty$.
\end{enumerate}
\end{definition}

\subsection{Function spaces for the overall error equation}

\begin{definition}
 $S^{l,\rho,m,\theta}$ is the set of all functions $f(x,Y)$ such that: \begin{enumerate}
\item $f$ is analytic inside $D(\rho) \times \Sigma(\theta, \frac{a}{\varepsilon})$ and  in $L^2 \left( \Gamma(\rho') \times \Gamma(\theta', \frac{a}{\varepsilon} \right)$ $\forall |\rho'|<\rho,$ $\forall |\theta'|<\theta$ ;
\item $|f|_{l,\rho,m,\theta} = \sum \limits_{k \leq m} |\partial_Y^k f|_{l-k,\rho,\theta} < \infty$.
\end{enumerate}
\end{definition}
\begin{definition} $S^{l,\rho,m,\theta}_{\beta,T}$ the set of all functions $f(t,x,Y)$ such that: \begin{enumerate}
\item $f \in B^0_\beta \left( [0,T], S^{l,\rho,m,\theta} \right)$;
\item $ |f|_{l,\rho,m,\theta,\beta,T}= \sup \limits_{t \in [0,T]} |f|_{l,\rho - \beta t, m,\theta - \beta t}< \infty$.
\end{enumerate}
\end{definition}
In the equation \eqref{edecompositione} of section \ref{sezionedellerrore}, we further decompose the error $\textbf{e}$: the space $S^{l,\rho,1,\theta}_{\beta,T}$ is the functional setting in which the ACK theorem can be applied for the term $\textbf{e}^*$, and the image of $\textbf{e}^*$ under the operator $\mathcal{N}^*$ defined in \ref{subsezioneNSoper} is in the space $L^{l,\rho,\theta}_{\beta,T}$, defined below.
\begin{definition}
$L^{l,\rho,\theta}_{\beta,T}$ is the set of functions $f(x,Y,t)$ such that

 $f \in C^0([0,t],S^{l,\rho-\beta t,0,\theta-\beta t})$,
  $\partial_Y f$, $\partial_{YY} f$, $\partial_t f$ $\in C^0([0,t],S^{l-2,\rho-\beta t,0,\theta-\beta t})$ $\forall t \leq T \leq \min \{\rho/\beta,0,\theta\}$, with norm \begin{equation}
\begin{split}
|f|_{l,\rho,\theta,\beta,T} = \sum \limits_{j=0}^{1} \sum \limits_{\alpha \leq l - 2 j} \sup \limits_{t \in [0,T]} |\partial_t^j \partial_x^\alpha f(\cdot,\cdot,t)|_{0,\rho-\beta t, \theta - \beta t} + \\ \sum \limits_{1 \leq \alpha_1 \leq 2} \sum \limits_{|\alpha_2| \leq l-2} \sup \limits_{0 \leq t \leq T} |\partial_Y^{\alpha_1} \partial_x^{\alpha_2} f(\cdot,\cdot,t)|_{0,\rho- \beta t, \theta - \beta t}.
\end{split} 
\end{equation}
\end{definition}
\subsection{Algebra properties and Cauchy estimates}
Let $d'$ be the dimension of $x$, i.e. $d'=1$ or $d'=2$ when solving NS equations in the half-plane or in the half-space respectively; moreover, call $d=d'+1$
\begin{lemma}\label{algebrax}
Let $u(x)$, $v(x)$ be in $H^{l,\rho}$, with $l > d'/2$. Then $uv \in H^{l,\rho}$, with \begin{equation}\label{algebraxformula}
|uv|_{l,\rho} \leq C |u|_{l,\rho} |v|_{l,\rho}.
\end{equation}
\end{lemma}
 \begin{lemma}\label{algebraHxy}
Let $u(x,y)$, $v(x,y)$ be in $H^{l,\rho,\theta}$, with $l > d/2$. Then $uv \in H^{l,\rho,\theta}$, with 
\begin{equation}\label{algebraHxyformula}
|uv|_{l,\rho,\theta} \leq C |u|_{l,\rho,\theta} |v|_{l,\rho,\theta}.
\end{equation}
\end{lemma}
One can easily verify both lemmas using the Paley-Wiener theorem and the usual argument used 
to prove the algebra properties of Sobolev spaces. 
%The lemmas can also be proven directly from the algebra properties of Sobolev spaces of functions of real variables, 
%obtaining the relations for fixed values of the imaginary parts of $x$ and $y$ and then taking the supremum: the main problem of this approach is that we take the supremum with respect to $Im(x)$ after the $L^2$ norm in $y$, but since an $L^p$ version of the three lines lemma holds \cite{BK2007}, we have 
%\begin{equation}
%\begin{split}
%\sup \limits_{|Im(x)|< \rho} \left||f|_{L^2(\Gamma(Im(x)))}\right|_{L^2(\Gamma(\theta',a))} \leq  
%\left| \sup \limits_{|Im(x)|< \rho} |f|_{L^2(\Gamma(Im(x)))}\right|_{L^2(\Gamma(\theta',a))} \leq \\  
%\left|  |f|_{L^2(\Gamma(\rho))}\right|_{L^2(\Gamma(\theta',a))} +  \left|  |f|_{L^2(\Gamma(-\rho))}\right|_{L^2(\Gamma(\theta',a))} \leq 
%2   \left|\sup \limits_{|Im(x)|< \rho}|f|_{L^2(\Gamma(Im(x)))}\right|_{L^2(\Gamma(\theta',a))}.
%\end{split}
%\end{equation}
%This approach shall also be used in the proof of lemma \ref{ydey}, and allows to obtain results in a general $L^p$ with $p \neq 2$. 

The use of complex variables allows to have simple estimates for the norms of the derivatives: using Cauchy formula for derivatives, we immediately have that \begin{equation}\label{Cauchyxest}
|\partial_x^k u|_{0,\rho'} \leq \frac{k!}{(\rho-\rho')^k} |u|_{0,\rho}.
\end{equation}
For functions holomorphic in a cone or in a conoid, we cannot bound the norm of the derivatives with the norm of the function in a larger cone with the same vertex, because we cannot use the Cauchy formula for derivatives with a fixed radius in all the path of integration. 
If we use a radius linearly growing with $|y|$ for $Re(y)<a$, constant for $Re(y)\geq a$, we easily obtain for a function holomorphic in $\Sigma \left( \theta,a \right)$ that,  $\forall p \in [1,+\infty]$ \begin{equation}\label{CauchyestyL2}
\sup \limits_{|\theta'|<\theta}|\min\{a,|y|\}\partial_yf(y)|_{L^p(\Gamma(\theta',a))} \leq  \frac{C}{\bar{\theta}-\theta} \sup \limits_{|\theta'|<\bar{\theta}} |f|_{L^p(\Gamma(\theta',a))}.
\end{equation}
From equation \eqref{CauchyestyL2}, we obtain the following lemma \cite{SC1998A}, which is crucial for the estimate of the nonlinear term in the Euler equations. 
\begin{lemma}\label{ydey}
Assume that $u$,$v$ $\in H^{l,\rho,\theta}$, with $u(y=0)=0$ and $l>d/2+1$. Then for any $\theta'<\theta$ we have \begin{equation}\label{ydeyeq}
|u \partial_y v|_{l,\rho,\theta'} \leq C |u|_{l,\rho,\theta'} \frac{|v|_{l,\rho,\theta}}{\theta-\theta'}.
\end{equation}
\end{lemma}
The Prandtl equations with initial data exponentially decaying in $Y$ need the following estimates, which can be found in \cite{SC1998A}.
\begin{lemma}
Let $f \in K^{l,\rho',\theta'',\mu''}$. Then \begin{equation}
|min\{1,|Y|\} \partial_Y f|_{l,\rho',\theta',\mu'} \leq \frac{|f|_{l,\rho',\theta'',\mu'}}{\theta''-\theta'} + \mu' |f|_{l,\rho',\theta',\mu'}
\end{equation}
\begin{equation}
|Y \partial_Y f|_{l,\rho',\theta',\mu'} \leq \frac{|f|_{l,\rho',\theta'',\mu'}}{\theta''-\theta'} + \mu' \frac{|f|_{l,\rho',\theta',\mu''}}{\mu''-\mu'} + |f|_{l,\rho',\theta',\mu'}
\end{equation}
\end{lemma}
Notice that a different exponential decay is needed in order to estimate the product of $\partial_Y f$ with a linearly increasing function. The Prandtl equations are still well posed when the initial data decay with $Y$ with a polynomial rate, see \cite{KV2013}, and also  \cite{CLS2013}; in this case, using a radius linearly growing with $|Y|$ in the Cauchy formula for derivatives, it is easy to obtain
\begin{equation}\label{CauchyestyLinfpesato}
\sup \limits_{Y \in \Sigma({\theta})}(1+|Y|)^\alpha |Y\partial_Yf(Y)| \leq  \frac{C}{\bar{\theta}-\theta} \sup \limits_{Y \in \Sigma({\bar{\theta}})} (1+|Y|)^\alpha |f(Y)|,
\end{equation}
which means that, unlike the exponential case, there is no need to change the polynomial rate of decay. 

\subsection{Paths of integration}\label{pathsubsec}
In the present paper, we shall solve equations involving heat operators in several instances: for example, boundary layer equations, where diffusion in the normal $Y$-variable appears,  or the error equation, where diffusion is both in the $x$- and $Y$-variable is present. 
The method we shall use is based on representing the solutions utilizing the convolution with the appropriate gaussian. 
As long as the integrand is holomorphic, one can deform the integration path, choosing the most convenient one.  
However, when estimating the solution and passing the modulus inside the integral, the integrand is no longer holomorphic, and the integral depends on the particular path of integration; 
therefore, in many estimates,  a careful a priori choice of the path of integration is helpful. 
Taking the convolution with respect to $x$, the best choice for the path is the one that makes the argument of the gaussian real, i.e., 
the path $x' \in \Gamma(Im(x))$.  
On the other hand, for the convolution in the $Y$ variable  between a gaussian and a function analytic in a conoid, given that the domain of analyticity shrinks near the origin, 
one cannot  take $Im(Y')=Im(Y)$; 
therefore, we choose the path consisting of the segment connecting the origin and $Y$, and the half-line from $Y$ to $+\infty + i Im(Y)$ parallel to the real axis. Said differently, one has: 
\begin{equation}\label{path}
\begin{split}
Y'=\left\{
\begin{array}{ll}
\frac{Y}{|Y|} r \, \, \, \, r \in [0,|Y|] \\
r+ i Y_{im} \, \, \, \, r \in [Y_r,+\infty[
\end{array}
\right. .
\end{split}
\end{equation}
In the first part of the above path of integration one has 
\begin{equation}\label{estpath1}
\left| e^{- \frac{(Y-Y')^2}{4(t-s)}} \right|= e^{- \frac{(|Y|-r)^2}{4(t-s)} \frac{Y_r^2-Y_{im}^2}{|Y|^2}} \leq e^{- \frac{(|Y|-r)^2}{4(t-s)} \cos(2 \theta)} ;
\end{equation}
to get  the last inequality we have used 
\begin{equation}
\frac{Y_r^2-Y_{im}^2}{Y_r^2+Y_{im}^2} = \frac{Y_r^2(1-\tan^2(\theta_Y))}{Y_r^2 (1 + \tan^2(\theta_Y))}=\cos(2 \theta_Y) \geq \cos(2 \theta),
\end{equation}
where $\theta_Y$ is the argument of the complex number $Y$. Analogously, one can show that 
\begin{equation}\label{estpath2}
\left| e^{- \frac{(Y+Y')^2}{4(t-s)}} \right| \leq e^{- \frac{(|Y|+r)^2}{4(t-s)} \cos(2 \theta)} ,
\end{equation}
\begin{equation}\label{estpath3}
\left| e^{- \frac{(Y')^2}{4s}} \right| \leq e^{- \frac{r^2}{4s} \cos(2 \theta)} ,
\end{equation}
where the latter inequality will be used for estimating the singular term $u^S$. 

In the second part of the path \eqref{path}, we have that  
\begin{equation}
\left| e^{- \frac{(Y-Y')^2}{4(t-s)}} \right|= e^{- \frac{(Y_r-r)^2}{4(t-s)} },
\end{equation}
while 
\begin{equation}
\left| e^{- \frac{(Y+Y')^2}{4(t-s)}} \right|= e^{ \frac{-(Y_r+r)^2+4 Y_{im}^2}{4(t-s)} } \leq e^{- (1-\tan^2(\theta)) \frac{(Y_r+r)^2}{4(t-s)} },
\end{equation}
where we have used 
$$
(2Y_{im})^2 \leq \tan^2 (\theta_Y) (2 Y_r)^2 \leq \tan^2(\theta) (Y_r+r)^2. 
$$
Finally, since $(Y_{im})^2 \leq \tan^2(\theta) r^2$, we have 
\begin{equation}
\left| e^{- \frac{(Y')^2}{4s}} \right|\leq e^{- (1-\tan^2(\theta)) \frac{r^2}{4s} }.
\end{equation}

Therefore, we conclude that after a change of variable, 
one can deduce the estimates involving the complex gaussian on the paths \eqref{path} by estimates involving the real gaussian; these estimates are up constants that blow up when $\theta \rightarrow \pi/4$.

For a function analytic in a cone, we take as path the half line which starts at the origin and passes through $Y$: therefore, the estimates \eqref{estpath1}, \eqref{estpath2} and \eqref{estpath3} are valid along all the path.

\section{The abstract Cauchy Kowalewski theorem}\label{sezioneteoremadicauchykow}

Consider the equation \begin{equation}\label{ACK equation}
u + F(t,u)=0 \, \, \, \, \, \, t \in [0,T].
\end{equation}
Let $\{ X_{\rho} \}_{\rho \in ]0,\rho_0]}$ a scale of Banach spaces, with $X_{\rho'} \subseteq X_{\rho''}$ and $|\cdot|_{\rho'} \leq |\cdot|_{\rho''}$ when $\rho'' \leq \rho' \leq \rho$
\begin{theorem}{(ACK theorem)}\label{ACK thm}
Suppose that $\exists R >0$, $\rho_0 >0$, $\beta_0 >0$ such that if $0<\tau \leq T \leq \rho_0 / \beta_0$ the following properties hold: \begin{enumerate}
\item $\forall 0 < \rho' < \rho \leq \rho_0 - \beta_0 \tau$ and $ \forall u$ such that $ \{ u \in X_\rho: \sup \limits_{t \in [0,\tau]} |u(t)|_\rho \leq R \}$ the map $F(u,t): [0,\tau] \mapsto X_{\rho'} $ is continuous;
\item $ \forall 0 < \rho \leq \rho_0- \beta_0 \tau$ the function $F(t,0):[0,\tau] \mapsto \{ u \in X_\rho : \sup \limits_{t \in [0,\tau]} |u(t)|_\rho \leq R \}$ is continuous and \begin{equation}
|F(t,0)|_\rho \leq R_0 < R;
\end{equation}
\item for any $\beta \leq \beta_0$, $0<\rho'<\rho(s) \leq \rho_0 - \beta_0 s$ and $u^1$ and $u^2 \in \{ u: u(t) \in X_{\rho_0 - \beta t }\leq R\}$ \begin{equation}
\begin{split}
|F(t,u^1)-F(t,u^2)|_{\rho'} \leq \\ C \left[ \int \limits_0^t  \frac{|u^1-u^2|_{\rho(s)}}{\rho(s)-\rho'} + \frac{|u^1-u^2|_{\rho'}}{(t-s)^{\alpha_1}} + \frac{1}{t^{\alpha^2}} \int \limits_0^t |u^1-u^2|_{\rho'} d s \right],
\end{split}
\end{equation}
with $\alpha_1<1$, $\alpha_2 <1$ and $C$ independent of $t, \tau, u^1, u^2, \beta, \rho', \rho(s)$.
\end{enumerate}
Then $\exists \beta > \beta_0$ such that equation \eqref{ACK equation} has a unique solution $u$ such that $\forall \rho \in ]0,\rho_0[$ $u(t) \in X_{\rho}$ $\forall t \in [0,(\rho_0-\rho)/\beta]$; moreover $ \sup \limits_{\rho + \beta t< \rho_0} |u(t)|_\rho \leq R$.
\end{theorem}
The spaces introduced in  \cref{sezionespazifunzionali} are Banach scales with respect to the parameters defining the complex domains and the exponential decay. The use of an analytic setting is mainly due to Prandtl's equations, since the well posedness results available for these equations use either an analytic setting \cite{KV2013}  and \cite{LCS2003} or some monotonicity assumption \cite{OS1997}. A key point of the ACK theorem is to prove that $F$ is a contraction in an auxiliary Banach space: this can be proved with slight modifications of the proof in \cite{LCS2003}.

\section{The main result}\label{sezionerisultatoprincipale}

We now state the main result of the paper. 
An informal statement was given by Theorem \ref{informalmain}, in section \ref{informal_main_result}, where 
the reader can find a detailed explanation of the meaning of the different  terms  
appearing in the asymptotic expansion \eqref{asymptotic_formal} below. 
\begin{theorem}\label{mainth}
Assume that $\textbf{u}_0 \in H^{l,\rho,\theta}$, with $\gamma_n \textbf{u}_0=0$, $l \geq 6$. 
Then, for any $\mu>0$, there exist $\bar{\rho}<\rho$, $\bar{\theta}<\theta$, $\bar{\beta}>0$, all independent of $\nu$, such that 
the solution of the Navier-Stokes equation \eqref{Navier-Stokes} can be written as
\begin{equation}
\textbf{u}=\textbf{u}^E + \bar{\textbf{u}}^P+ \varepsilon [ \textbf{u}^E_{(1)}+\bar{\textbf{u}}^P_{(1)}+ \textbf{e}], \label{asymptotic_formal}
\end{equation}
where: 
\begin{enumerate}
\item The term $\textbf{u}^E \in H^{l,\bar{\rho},\bar{\theta}}_{\bar{\beta},T}$ is the solution of the Euler equations \eqref{Euler}.
\item The term  $\bar{\textbf{u}}^P$ is the modified Prandtl solution given by \cref{utildeP} and \cref{vbarP}.  
The following decomposition in regular and singular part holds:
$$
\bar{\textbf{u}}^P=\bar{\textbf{u}}^R+ \bar{ \textbf{u}}^S, 
$$
where the regular term $\bar{\textbf{u}}^R \in K^{l,\bar{\rho},\bar{\theta},\mu}_{\bar{\beta},T}$, and the singular term $\textbf{u}^S$ is given by 
\eqref{u_singular} and \eqref{v_singular}. 
\item The term  $\textbf{u}^E_{(1)}$ is the first order correction to the inviscid flow solving the system \eqref{CorrEul}.
The following decomposition in regular and singular part holds:
 $$
 \textbf{u}^E_{(1)}=\textbf{u}^{ER}_{(1)}+\sqrt{t} \boldsymbol{w}^S_b(x,y), 
 $$
where $\textbf{u}^{ER}_{(1)} \in  \tilde{H}^{l,\bar{\rho},\bar{\theta}}_{\bar{\beta},T,1}$ and 
$\boldsymbol{w}^S_b \in H^{l,\bar{\rho},\bar{\theta}}$ is given by \eqref{wboundary}.
\item $\bar{\textbf{u}}^P_{(1)}=\left(\bar{u}^P_{(1)}, \varepsilon \bar{v}^P_{(1)}\right) $ is the first order correction to the boundary layer flow; 
 $\bar{u}^P_{(1)}$ solves the linear heat equation \eqref{uP1syst}, while $ \bar{v}^P_{(1)}$ is given by the incompressibility condition \eqref{vP1}. 
The following decomposition in regular and singular part holds:
\begin{equation}
\bar{\textbf{u}}^P_{(1)}= \bar{\textbf{u}}^R_{(1)} + \bar{\textbf{u}}^S_{(1)} \label{deco_prandtl_first_order}
\end{equation}
where the regular part $\bar{\textbf{u}}^R_{(1)} \in K^{l,\bar{\rho},\bar{\theta},\mu}_{\bar{\beta},T}$ 
and the singular part $\bar{\textbf{u}}^S_{(1)}$ is given by \eqref{us1syst}. 
\item The term $\textbf{e}$ is an overall error that closes the asymptotic procedure with $\textbf{e} \in L^{l,\bar{\rho},\bar{\theta}}_{\bar{\beta},T}$.
\end{enumerate}
\end{theorem}
The result concerning the Euler flow can be found in \cite{SC1998A} as Theorem 4.1. 
The result for $\bar{\textbf{u}}^P$ can be easily obtained with a slight modification of the argument used in \cite{CLS2013}. 

\begin{remark}
	One can give a more detailed description of the structure of $\bar{\textbf{u}}^P$ as follows. 
The solution $u^P$ of Prandtl equations with non compatible data can be written as the sum of 
\begin{enumerate}
\item a singular term $u^S$ given by \eqref{u_singular};
\item the solution $u^D$ of Prandtl's equations with compatible data ($\gamma u^D= u_{00}$) and an initial datum given by the initial value of the Euler flow at the boundary, which  means that $u_D$ approaches the Euler flow for $Y \rightarrow + \infty$ at any time with any exponential rate;
\item an interaction term, which is exponentially decaying: indeed, the interaction part can be written in terms of the heat operator $E_2$ for the half space with zero boundary and initial condition and a forcing term, so in $Y$ we have a convolution between a gaussian and some terms deriving from $\tilde{u}^D$ and $u^S$. Now we have 
\begin{equation}
e^{\mu Re(Y)}= e^{\mu Re(Y-Y')} e^{\mu Re(Y')},
\end{equation}
with $e^{\mu Re(Y')}\tilde{u}^D$ bounded, while \begin{equation}
\begin{split}
e^{\mu Re(Y-Y')} \left| \frac{e^{-\frac{(Y-Y')^2}{4(t-s)}}}{\sqrt{4\pi (t-s)}} \right| \leq e^{\mu Re(Y-Y')} e^{-c\frac{Re(Y-Y')^2}{8T}} \frac{e^{-\frac{cRe(Y-Y')^2}{8(t-s)}}}{\sqrt{4\pi (t-s)}} \leq \\ \leq C \frac{e^{-\frac{cRe(Y-Y')^2}{8(t-s)}}}{\sqrt{4\pi (t-s)}},
\end{split}
\end{equation}
with $c$ depending on $\theta$. A similar argument applies to $e^{\mu Re(Y)}u^S$.
\end{enumerate} 
\end{remark}
The result for $\bar{\textbf{u}}^P$ can be obtained with minor modifications of the argument in \cite{CLS2013}, taking into account that the decay in $Y$ is exponential (rather than polynomial) and that we have analyticity also in the normal variable.
\begin{remark}
To prove  \cref{mainth}, we use the fact that when we reduce the strip of $x$-analyticity, we can get as much $x$-regularity 
 as as needed. 
 However, while we can obtain arbitrary regularity in the tangential variable, this does not hold for the normal variable 
 since the reduction of the cone of analyticity does not provide any additional regularity: 
 the regularity in the normal variable does not exceed the regularity of the initial datum.
\end{remark}

The construction of $\textbf{u}^E$ and $\bar{\textbf{u}}^P$ already appears in the literature, see \cite{SC1998A} and \cite{CLS2013} respectively; 
therefore, for the rest of the paper, we shall construct the remaining terms appearing in the asymptotic expansion \eqref{asymptotic_formal}.

\section{Correction to the Euler flow}\label{sezionecorrezioneeulero}
The first-order correction to the Euler flow is produced by the inflow at the boundary generated by the Prandtl solution; 
given that the Prandtl solution has a regular and a singular part, in \eqref{VPrandtlboundary}
we have split the inflow 
$g$ in singular and regular part, $g^S$ and $g^R$, respectively.  
The first order Euler equations \eqref{CorrEul} are linear; therefore, by the superposition principle,  
we can decompose the first order Euler correction $\bar{\textbf{u}}^E_{(1)}$ as in
\eqref{euler_first_decomposition}, i.e., as singular and regular part. 
The singular part $\boldsymbol{w}^S$ solves
\begin{equation}\label{CorrEulSing}
\begin{split}
\partial_t \boldsymbol{w}^S + \boldsymbol{w}^S \cdot \nabla \textbf{u}^E + \textbf{u}^E \cdot \nabla \boldsymbol{w}^S + \nabla p^{w^S} = \textbf{0}, \\
\nabla \cdot \boldsymbol{w}^S=0, \\
\gamma_n \boldsymbol{w}^S= g^S=-\sqrt{\frac{4}{\pi}} \sqrt{t} \partial_x u_{00}(x) ,\\
\boldsymbol{w}^S\big|_{t=0}= \textbf{0},
\end{split}
\end{equation}
while the regular part $\boldsymbol{w}^R$ solves  the system obtained by \eqref{CorrEulSing} replacing $g^S$ with $g^R$. 
The regularity of $g^R$ is the same that the boundary value of the normal part of the Prandtl's correction would have in the compatible case; 
therefore, we have $\boldsymbol{w}^R \in  \tilde{H}^{l-1,\bar{\rho},\bar{\theta}}_{\bar{\beta},T,1}$. 
\begin{remark}
The regularity in the $y$-variable is greater than the one stated in \cite{SC1998B}; however, for this term, 
the regularity in the normal variable depends only on "symmetric" operators, which allow a number of derivatives in $y$ equal to the ones that have not been used in $x$; 
this enhanced regularity will allow proving that some terms which are zero at time $t=0$ are $O(t)$, 
which will turn out to be helpful for the estimates of some nonlinear terms in the equation of $\textbf{e}$.\end{remark}

We now pass to the analysis of the singular part $\boldsymbol{w}^S$. 
If one denotes by $\xi'$ the Fourier variable corresponding to $x$,  it is useful to write the solution as 
\begin{equation}\label{omegaessestar}
\boldsymbol{w}^S(t,\xi',y)=\sqrt{t} \boldsymbol{w}^S_b   + I_{\sqrt{t}} \boldsymbol{w^{S^*}} ,
\end{equation}
where
\begin{equation}\label{wboundary}
\boldsymbol{w}^S_b(\xi',y)=C_1 (- i \frac{\xi'}{|\xi'|},1) e^{-|\xi'|y} i \xi' u_{00}(\xi'),
\end{equation}
is divergence-free and takes into account the normal inflow. 
The operator $I_{\sqrt{t}}$ is a weighted integration in time: 
 \begin{equation}
I_{\sqrt{t}} \boldsymbol{w^{S^*}} = \int \limits_0^t \sqrt{\tau} \boldsymbol{w}^{S^*}(\tau,\xi',y) d \tau.
\end{equation}
The rest of this section is devoted to the construction $\boldsymbol{w^{S^*}}$. Call $P$ the Leray projector 
(projection on the space of divergence-free $L^2$ functions with zero normal component at the boundary); 
 substituting the ansatz \eqref{omegaessestar} into \eqref{CorrEulSing}, and then applying the Leray projector, 
 one derives the following equation for $\boldsymbol{w^{S^*}}$,
\begin{equation}\label{omegaSstar} \begin{split}
\boldsymbol{w}^{S^*}+P  \left[ \frac{I_{\sqrt{t}} \boldsymbol{w^{S^*}}}{\sqrt{t}} \cdot \nabla \textbf{u}^E + \textbf{u}^E \cdot \nabla\frac{I_{\sqrt{t}} \boldsymbol{w^{S^*}}}{\sqrt{t}} + \boldsymbol{w}^S_b \cdot \nabla \textbf{u}^E  +   \textbf{u}^E \cdot \nabla \boldsymbol{w}^S_b \right] = \textbf{0},
\end{split}
\end{equation}
which is a suitable form for the application of the abstract Cauchy-Kowalewski theorem. 
Since $P$ is a bounded operator in $H^{l,\rho,\theta}$ (see \cite{SC1998A}) and commutes with time derivatives, we can estimate $\textbf{v}$ instead of $P \textbf{v}$ for any $\textbf{v}$. We shall show that, assuming that $u_{00} \in H^{l+1,\rho}$ and $\textbf{u}^{E} \in H^{l,\rho,\theta}_{\beta,T}$, then $\boldsymbol{w}^{S*} \in \tilde{H}^{l-1,\rho,\theta}_{\bar{\beta},T,1}$ for some $\bar{\beta}>\beta$: such a regularity for $u_{00}$ can always be obtained from the trace of a function $u_0 \in H^{l,\rho,\theta}$ by reducing the strip of analyticity.

\subsection{The forcing term}
For any $\rho'<\rho-\beta t$, $\theta' \leq \theta-\beta t$, we use the algebra property \eqref{algebraHxyformula} to obtain \begin{equation}
|\boldsymbol{w}^S_b \cdot \nabla \textbf{u}^E|_{l-1,\rho',\theta'} \leq C |\boldsymbol{w}^S_b|_{l-1,\rho',\theta'} |\textbf{u}^E|_{l,\rho',\theta'} \leq C |u_{00}|_{l,\rho} |\textbf{u}^E|_{l,\rho,\theta,\beta,T}
\end{equation} 
and similarly \begin{equation}
|\textbf{u}^E \cdot \nabla \boldsymbol{w}^S_b |_{l-1,\rho',\theta'} \leq C |u_{00}|_{l+1,\rho} |\textbf{u}^E|_{l-1,\rho,\theta,\beta,T}.
\end{equation}
\subsection{Quasi contractiveness}
For the term $\textbf{u}^E \cdot \nabla\frac{I_{\sqrt{t}} \boldsymbol{w^{S^*}}}{\sqrt{t}}$, we have \begin{equation}
 \begin{split}
\left|u^E \partial_x \frac{I_{\sqrt{t}} \boldsymbol{w}^{S*}}{\sqrt{t}}\right|_{l-1,\rho',\theta'} \leq C |u^E|_{l-1,\rho',\theta'} \int \limits_0^t \frac{\sqrt{\tau}}{\sqrt{t}} |\partial_x \boldsymbol{w}^{S^*}|_{l-1,\rho',\theta'} d \tau \leq\\ \leq  C |u^E|_{l-1,\rho,\theta,\beta,T} \int \limits_0^t \frac{|\boldsymbol{w}^{S*}|_{l-1,\rho(s),\theta'}}{\rho(s)-\rho'} ds,
\end{split}
\end{equation} while, with an argument similar to the one used in the proof of lemma \ref{ydey}, since $ \partial_y \frac{I_{\sqrt{t}} \boldsymbol{w^{S^*}}}{\sqrt{t}}$ multiplies $v^E$, which goes to zero linearly as $y$ approaches zero, we have \begin{equation}
\left| v^E \partial_y \frac{I_{\sqrt{t}} w^{S^*}}{\sqrt{t}} \right|_{l-1,\rho',\theta'} \leq C |v^E|_{l-1,\rho,\theta,\beta,T} \int \limits_0^t \frac{|\boldsymbol{w}^{S^*}|_{l-1,\rho',\theta(s)}}{\theta(s)-\theta'}ds.
\end{equation}
The term $\frac{I_{\sqrt{t}} \boldsymbol{w}^{S^*}}{\sqrt{t}} \cdot \nabla \textbf{u}^E$ is easier to estimate, since no Cauchy estimate is needed \begin{equation}
|\frac{I_{\sqrt{t}} \boldsymbol{w}^{S^*}}{\sqrt{t}} \cdot \nabla \textbf{u}^E| \leq C |\textbf{u}^E|_{l,\rho,\theta,\beta,T} \int \limits_0^t |\boldsymbol{w}^{S*}|_{l-1,\rho',\theta'} d s.
\end{equation} 
By the ACK theorem, we obtain the existence of a unique solution of equation \eqref{omegaSstar} in $\tilde{H}^{l-1,\rho,\theta}_{\bar{\beta},T,0}$; taking the time derivative of that equation, it's easy to see that $\boldsymbol{w}^{S*} \in \tilde{H}^{l-1,\rho,\theta}_{\bar{\beta},T,1}$.
\section{Boundary layer corrector}\label{sezionecorrezioneprand}
The first order correction $\bar{\textbf{u}}^P_{(1)}=(u^P_{(1)}, \varepsilon \bar{v}^P_{(1)})$ is the sum (see \eqref{deco_prandtl_first_order}) 
of the solution $\bar{\textbf{u}}^S_{(1)}$ of the system \begin{equation}\label{us1syst}
\begin{split}
(\partial_t - \partial_{YY}) u^S_{(1)}= u^S, \\
\gamma u^S_{(1)}=-\gamma \boldsymbol{w}^S, \\
u^S_{(1)}(t=0)=0, \\
\bar{v}^S_{(1)}= \varepsilon \int \limits_Y^{+\infty} \partial_x u^S d Y',
\end{split}
\end{equation}
and the solution $\bar{\textbf{u}}^R_{(1)}$ of the system obtained by \eqref{us1syst} replacing $\boldsymbol{w}^S$ with $\boldsymbol{w}^R$
We have $\bar{\textbf{u}}^R_{(1)} \in K^{l,\bar{\rho},\bar{\theta},\mu}_{\bar{\beta},T}$, see \cite{SC1998B} for proof. In Fourier terms, $u^S_{(1)}$ is given by
\begin{equation}
u^S_{(1)}=- C_1 |\xi'| u_{00}(\xi') F_{1/2}(t,Y)
\end{equation} 
where \begin{equation}
F_{1/2}= \int \limits_0^t \frac{Y}{2(t-s)} \frac{e^{- \frac{Y^2}{4(t-s)}}}{\sqrt{4 \pi (t-s)}} \sqrt{s} ds= \frac{1}{\sqrt{\pi}} \int \limits_{\frac{Y}{2\sqrt{t}}}^{+\infty} e^{- \sigma^2} \sqrt{t- \frac{Y^2}{4 \sigma^2}} d \sigma,
\end{equation}
We estimate $u^S_{(1)}$ in terms of $F$: we have \begin{equation}
|F_{1/2}| \leq C e^{- \frac{Y^2}{8t}} \sqrt{t}.
\end{equation}
The derivative with respect to $Y$ is \begin{equation}
\partial_Y F_{1/2}=- \frac{1}{\sqrt{\pi}} \int \limits_{\frac{Y}{2\sqrt{t}}}^{+ \infty} e^{- \sigma^2} \frac{Y}{4 \sigma^2} \frac{d \sigma}{\sqrt{t- \frac{Y^2}{4 \sigma^2}}}= - \frac{2}{\sqrt{\pi}} \int \limits_0^t e^{- \frac{Y^2}{4(t-s)}} \frac{ds}{\sqrt{s} \sqrt{(t-s)}},
\end{equation}
so \begin{equation}
| \partial_Y F_{1/2}| \leq C e^{- \frac{Y^2}{4t}}.
\end{equation}
The term $\partial_{YY} F_{1/2}$ is singular near the boundary (although $Y \partial_{YY} F_{1/2} $ is not).
\\ The normal component $\bar{v}^S_{(1)}$ has a better regularity with respect to $t$ and $Y$; we have \begin{equation}
|\partial_Y^j \bar{v}^S_{(1)}|_{k,\rho} \leq C |u_{00}|_{k+2,\rho} e^{- c\frac{|Y|^2}{t}} t^{1-\frac{j}{2}}
\end{equation}

\section{The remainder $\textbf{e}$}\label{sezionedellerrore}

The remainder $\textbf{e}$ satisfies equations  \eqref{corrector_e}, which are Navier-Stokes type equations with source term and 
non-homogeneous boundary conditions. 
It is important to notice that the source term contains $\partial_{xx}u^S$ which is the most singular term.
Therefore, we decompose the error $\textbf{e}$ as 
\begin{equation}\label{edecompositione}
\textbf{e}= \mathcal{N}^*\textbf{e}^* + \boldsymbol{\sigma} + \textbf{h}.
\end{equation}
In the above decomposition, the tangential part of $\textbf{h}=(h',\varepsilon h_n)$ takes care of the most singular term $\partial_{xx}u^S$; in fact, 
$h'$ solves the heat equation in the half space with forcing term $\varepsilon \partial_{xx} u^S$, and 
homogeneous boundary and initial conditions,  see the system \eqref{eq::h_t}; 
 the normal component $h_n$, is obtained imposing the incompressibility condition and the decay at infinity in $Y$, 
 see \eqref{expression_h_n}. Therefore $h_n$ has a non zero trace at the boundary. 
 The construction of $\textbf{h}$ will be accomplished in subsection \ref{sub::heat_term}, while the fact that $\textbf{h} \in L^{l,\bar{\rho},\theta}_{\beta,T}$ 
 is stated in Proposition \ref{hinL}. 
 
 The term $\boldsymbol{\sigma}$ solves the Stokes equations with boundary conditions,  
 see the system \eqref{sigmasing} below. 
We have introduced the term $\boldsymbol{\sigma}$ to take into account  
 the boundary conditions of $\textbf{e}$, deriving from the boundary layer corrector, and the boundary conditions generated by $\textbf{h}$.
 The construction of $\boldsymbol{\sigma}$ will be accomplished in subsection \ref{sub::sigma}, while the fact that 
 $\boldsymbol{\sigma} \in L^{l,\bar{\rho},\theta}_{\beta,T}$ 
 is stated in Proposition \ref{sigmainL}. 
% We also remark that the 
% eventual initial conditions for $\textbf{e}$ (that in our case are homogeneous),

 Finally, $\mathcal{N}^*\textbf{e}^*$ solves the Stokes equations with homogeneous boundary conditions and with forcing term $\textbf{e}^*$, 
 see the system \eqref{Stokes_source}.    
 The Navier-Stokes operator $\mathcal{N}^*$ can be defined explicitly in terms of the Leray projector and and heat operator, 
 see the formula \eqref{NSoperator}. 
 The construction of the operator $\mathcal{N}^*$, and the necessary estimates, will be presented in subsection \ref{subsezioneNSoper}

To construct $\textbf{e}^*$ we shall use the ACK theorem in the $S^{l,\rho,1,\theta}$ setting: 
given that  $u^S \notin S^{l,\rho,1,\theta}$, this has led us to isolate 
the effect of $\varepsilon \partial_{xx} u^S$ introducing the term $\textbf{h}$.
The construction of $\textbf{e}^*$ will be accomplished in subsection \ref{sub::error_star}. 

 For the compatible case, in \cite{SC1998B}, the authors showed that $\textbf{e} \in L^{l,\rho,\theta}_{\beta,T}$ 
 by proving that both $\boldsymbol{\sigma}$ and $\textbf{e}^*$ were in that space.  
 In our case, $\boldsymbol{\sigma}$ is still in $ L^{l,\rho,\theta}_{\beta,T}$; however, we can say that $\textbf{e}^* \in S^{l,\rho,1,\theta}$, only; 
 nevertheless, we shall prove that the image under $\mathcal{N}^*$ is still in $L^{l,\rho,\theta}_{\beta,T}$.

\subsection{Heat term}  \label{sub::heat_term}
The tangential part $h'$ of $\textbf{h}$ satisfies 
\begin{equation} \label{eq::h_t}
\begin{split}
(\partial_t  - \partial_{YY})h' = \varepsilon \partial_{xx} u^S, \\
\gamma h'=0, \\
h'(t=0)=0.
\end{split}
\end{equation}
One can write the explicit expression of $h'$:
\begin{equation}
h'= \varepsilon \partial_{xx} u_{00} F(t,Y),  \label{expression_h_prime}
\end{equation}
where
\begin{equation}
F= \frac{2}{\sqrt{\pi}} \int \limits_0^t ds \int \limits_0^{+\infty} (E_0^--E_0^+) \int \limits_{\frac{Y'}{2\sqrt{s}}}^{+\infty} e^{-\sigma^2} d \sigma d Y',
\end{equation}
and $E_0^-$ and $E_0^+$ are:
\begin{equation}
E_0^-= \frac{e^{-\frac{(Y-Y')^2}{4(t-s)}}}{\sqrt{4 \pi (t-s)}}; \, \, E_0^+= \frac{e^{-\frac{(Y+Y')^2}{4(t-s)}}}{\sqrt{4 \pi (t-s)}}.
\end{equation}
The normal part of $\textbf{h}$ is obtained through the incompressibility condition expressed in the $Y$ variable; thus, it is a $O(\varepsilon)$ with respect to $h'$. We shall denote the normal part of $\textbf{h}$ as $\varepsilon h_n$ to stress this fact. We have 
\begin{equation}
	h_n= \int \limits_Y^{+\infty} \partial_{x}h' dY'.   \label{expression_h_n}
\end{equation}
One can give the following estimate:
\begin{proposition} \label{prop:estimate_h} For $j=0,1,2$, $\forall \mu >0$, we have
\begin{equation}
\begin{split}
\sup \limits_{Y \in \Sigma{(\theta)}} e^{\mu Re(Y)} | \partial_{Y}^j h' |_{l-2-j,\rho} \leq C \varepsilon t^{\frac{2-j}{2}} |u_{00}|_{l-j,\rho}, \\
\sup \limits_{Y \in \Sigma{(\theta)}} e^{\mu Re(Y)} | \partial_{Y}^j h_n |_{l-3-j,\rho} \leq C \varepsilon t^{\frac{3-j}{2}} |u_{00}|_{l-j,\rho}.
\end{split}
\end{equation} 
\end{proposition} 
The proof of the proposition is in appendix \ref{proof_on_h}. 
The exponential decay in the $Y$ variable stated in proposition \ref{prop:estimate_h} implies the boundedness in $L^2_Y$; 
moreover  every $u_0 \in H^{l,\rho,\theta}$ has trace $u_{00} \in H^{l-1,\rho} \subset H^{k,\bar{\rho}}$ $\forall k \geq l$, $\forall \bar{\rho} < \rho$. 
Therefore, we obtain the following result: 
%n.b. : u_0 è volutamente non in grassetto, in quanto è solo la parte tangenziale che stiamo considerando, essendo l'unica a dare origine a u_00
\begin{proposition}\label{hinL}
Assume $u_0 \in H^{l,\rho,\theta}$: then $\forall \bar{\rho}<\rho$, $\forall T$, $\forall \beta$, we have $\textbf{h} \in L^{l,\bar{\rho},\theta}_{\beta,T}$.
\end{proposition}

The estimates given in Proposition \ref{prop:estimate_h}  imply an $L^{\infty}_Y$ boundedness, 
which will be useful in the treatment of the nonlinear terms in the equation of $\textbf{e}^*$.

\subsection{Boundary value of the error} \label{sub::sigma}
The term $\boldsymbol{\sigma}=(\sigma_1, \sigma_2)$ is needed to cancel the boundary value of $h_n$ and $\bar{v}^P_{(1)}$; 
therefore, it satisfies 
\begin{equation}\label{sigmasing}
\begin{split}
(\partial_t-\partial_{YY}) \boldsymbol{\sigma} + \nabla \phi=0, \\
\nabla \cdot \boldsymbol{\sigma}=0,, \\
\gamma \boldsymbol{\sigma}=(0,\varepsilon G), \\
\boldsymbol{\sigma}(t=0)=\textbf{0},
\end{split}
\end{equation}
where \begin{equation}\label{Gespressione}
G= - \gamma \bar{v}^R_{(1)} -\gamma \bar{v}^P_{(1)} -\gamma h_n=  \partial_x \int_0^{+\infty} \left(u^R_{(1)}+ u^P_{(1)}+ h' \right)dY=\partial_x \tilde{G}.
\end{equation}
The above system is a Stokes problem with boundary datum and homogeneous initial datum. 
In \cite{SC1998B} one can find the procedure to solve such a problem; the solution writes as 
\begin{equation}\label{sigmaespressione}
\begin{split}
\sigma_1= \varepsilon N' E_1 G - N' \varepsilon \int \limits_0^Y \varepsilon |\xi'| e^{-\varepsilon |\xi'|(Y-Y')} E_1 G d Y' - N' e^{-\varepsilon |\xi'|Y} G, \\
\sigma_2 = \varepsilon e^{- |\xi'| \varepsilon Y} G + \varepsilon \int \limits_0^Y \varepsilon |\xi'| e^{-\varepsilon |\xi'|(Y-Y')} E_1 G,
\end{split}
\end{equation}
where 
\begin{equation}
N'= \frac{i \xi'}{|\xi'|}, \label{Riesz_Np}
\end{equation}
while $E_1f(t)$ gives the solution of the heat equation with boundary datum $f(t)$ and homogeneous initial datum; the explicit expression is
\begin{equation}
E_1 f = \int \limits_0^t \frac{Y}{2(t-s)} \frac{e^{-\frac{Y^2}{4(t-s)}}}{\sqrt{4 \pi (t-s)}} f(s) ds.
\end{equation}
\begin{proposition}\label{sigmainL}
	Let $\boldsymbol{\sigma} =\left(\sigma_1,\sigma_2\right)$ with $\sigma_i$ given in \cref{sigmaespressione}. 
Then $\boldsymbol{\sigma} \in L^{l,\bar{\rho}, \bar{\theta}}_{\bar{\beta},\bar{T}}$ and,  $\forall |\theta'| <\theta, \, \forall \theta \in ]0,\pi/4[, \, \forall a$
the following estimates, hold: 
 \begin{equation}\label{stimesigmaL2}
\begin{split}
||\boldsymbol{\sigma}|_{k,\rho}|_{L^2(\Gamma(\theta',a))} \leq C|\tilde{G}|_{k+1,\rho} , \\
||\partial_Y \boldsymbol{\sigma}|_{k,\rho}|_{L^2(\Gamma(\theta',a))} \leq C \varepsilon \left[ t^{3/4}|\partial_t G|_{k,\rho} + | G|_{\rho,k+1} \right], \\
||\partial_{YY} \boldsymbol{\sigma}|_{k,\rho}|_{L^2(\Gamma(\theta',a))} \leq C \varepsilon \left[ |G|_{k+2,\rho} +  t^{1/4}| \partial_t G|_{k+1,\rho} \right], \\
||\partial_t \boldsymbol{\sigma}|_{k,\rho}|_{L^2(\Gamma(\theta',a))} \leq C   \left| \partial_t \tilde{G} \right|_{k+1,\rho} ,
\end{split}
\end{equation}
\begin{equation}\label{stimesigmaLinf}
\begin{split}
\sup \limits_{Y} |\boldsymbol{\sigma}|_{k,\rho} \leq C |G|_{k,\rho}, \\
\sup \limits_{Y}| \partial_Y \boldsymbol{\sigma}|_{k,\rho} \leq C\varepsilon \left[ |G|_{k+1,\rho} + t^{1/2} |\partial_t G|_{k,\rho} \right], \\ 
\sup \limits_{Y}| \partial_{YY} \boldsymbol{\sigma}|_{k,\rho} \leq C\varepsilon \left[  |G|_{k+2,\rho}  |\partial_t G|_{k+1,\rho} \right].
\end{split}
\end{equation}
\end{proposition}
The proof of the proposition is straightforward, and can be obtained by applying Young's convolution inequality to the expressions \eqref{sigmaespressione}.
\begin{remark}
In the expression of $\boldsymbol{\sigma}$, 
the term $N'e^{-\varepsilon |\xi'|Y} G$ is the only one which is not $O(\varepsilon)$; 
however,  $\partial_Y \boldsymbol{\sigma}$ is $O(\epsilon)$, which implies that $\partial_y \boldsymbol{\sigma}$ is $O(1)$. 
\end{remark}
\begin{remark}
The estimates for $\boldsymbol{\sigma}$ and $\partial_Y \boldsymbol{\sigma}$ show that, near $t=0$, they are   small in $t$; 
in fact, $G(t=0)=0$ and continuosly differentiable in $t$, so that $G$ goes linearly to zero with $t$. 
The regularity of $G$ with respect to time is due to the regularizing effect that the integration in $Y$, appearing in \eqref{Gespressione}, 
has on "gaussian"-type singularities. 
\end{remark}
These remarks and the estimates in the $L^\infty_Y$ norm in \eqref{stimesigmaLinf} will be useful in the estimates for the nonlinear terms 
in the equation of $\textbf{e}^*$.

\subsection{The Navier-Stokes operator}\label{subsezioneNSoper}

In the present subsection we shall introduce the operator $\mathcal{N}^*$ and give the estimate Proposition \ref{NstarvinL}, 
which is the main result of the subsection. 
The Navier-Stokes operator $\mathcal{N}^*$ solves the time-dependent Stokes equation with a forcing term: 
$\boldsymbol{w}=\mathcal{N}^* \boldsymbol{w}^*$ is the solution of 
\begin{equation}\label{Stokes_source}
\begin{split}
(\partial_t-\varepsilon^2 \partial_{xx}- \partial_{YY}) \boldsymbol{w} + \nabla p^{w}= \boldsymbol{w}^* \\
\nabla \cdot \boldsymbol{w}= 0, \\
\gamma \boldsymbol{w}=\textbf{0}, \\
\boldsymbol{w}(t=0)=\textbf{0}.
\end{split}
\end{equation}
One can write the explicit expression of $\mathcal{N}^*$ in terms of a projection operator $ \bar{P}^\infty$, of the heat operator $ \tilde{E}_2$, 
and of the operator $ \mathcal{S}$ which solves the Stokes equations with boundary data:
\begin{equation}\label{NSoperator}
\mathcal{N}^*= \bar{P}^\infty \tilde{E}_2 - \mathcal{S} \gamma \bar{P}^{\infty} \tilde{E}_2. 
\end{equation}
In the  rest of the section we shall give the explicit expression of $\bar{P}^\infty$,  $\tilde{E}_2$, and   $\mathcal{S}$, 
and state the necessary estimates for $\mathcal{N}^*$. 

If $\textbf{v}$ is a vector field defined on the upper plane, $\bar{P}^\infty \textbf{v}$ is obtained in the following way: extend $\textbf{v}$ oddly for $Y<0$, then apply the Leray projector for functions defined on the whole space, and finally restrict the result to $Y \geq 0$. The explicit expression of the normal component is given by \begin{equation}
\begin{split}
\bar{P}^{\infty}_n \textbf{v} = \frac{1}{2} \varepsilon |\xi'| \Bigg[ \int \limits_0^Y dY' \left(e^{-\varepsilon|\xi'|(Y-Y')}-e^{-\varepsilon|\xi'|(Y+Y')} \right) (- N' v_1 + v_2) + \\
+ \int \limits_Y^{+\infty} d Y' \left( e^{\varepsilon|\xi'|(Y-Y')}(N' v_1+v_2) - e^{-\varepsilon|\xi'|(Y+Y')} (-N' v_1+ v_2) \right) \Bigg],
\end{split}
\end{equation} 
where the Riesz-type operator $N'$ is defined in \eqref{Riesz_Np}. 
The tangential component is given by 
\begin{equation}
\begin{split}
\bar{P}^{\infty '} \textbf{v}= v_1 - \frac{\varepsilon|\xi'|}{2} \Bigg[ \int \limits_0^Y d Y' \left( e^{-\varepsilon|\xi'|(Y-Y')}- e^{- \varepsilon |\xi'| (Y+Y')} \right) (v_1+N' v_2) + \\ + \int \limits_Y^{+\infty} d Y' \left( e^{\varepsilon|\xi'|(Y-Y')} (v_1-N'v_2)-e^{-\varepsilon|\xi'|(Y+Y')} (v_1+N' v_2) \right) \Bigg].
\end{split}
\end{equation}
Notice that, if we had extended the tangential part evenly and the normal part oddly, we would have obtained the Leray projector for the half space.

The operator $\tilde{E}_2$ is such that $\tilde{E}_2f$ solves the heat equation in the half space with source $f$ and homogeneous data:
\begin{equation}
\begin{split}
(\partial_t - \varepsilon^2 \partial_{xx} - \partial_{YY})u=f \\
\gamma u=0\\
u(t=0)=0
\end{split}
\end{equation}
The explicit expression can be  given in terms of convolutions with gaussians: 
\begin{equation} \begin{split}
\tilde{E}_2 f (x,Y,t)= \int \limits_0^t ds \int \limits_{-\infty}^{+\infty} \frac{e^{-\frac{(x-x')^2}{4(t-s)\varepsilon^2}}}{\sqrt{4 \pi \varepsilon^2 (t-s)}} d x' \\ \int \limits_0^{+\infty} \left( \frac{e^{-\frac{(Y-Y')^2}{4(t-s)}}}{\sqrt{4 \pi (t-s)}} - \frac{e^{-\frac{(Y+Y')^2}{4(t-s)}}}{\sqrt{4 \pi (t-s)}} \right)  f(x',Y',s) d Y'.
\end{split}
\end{equation}

The Stokes operator $\mathcal{S}$ is the operator such that $\mathcal{S} \textbf{g}$ solves the Stokes equations with 
boundary datum $ \textbf{g}$: 
\begin{equation}
\begin{split}
(\partial_y - \varepsilon^2 \Delta) \textbf{u}^S + \nabla p^S=0, \\
\nabla \cdot \textbf{u}^S=0, \\
\gamma \textbf{u}^S=\textbf{g}, \\
\textbf{u}^S(t=0)= \textbf{0}.
\end{split}
\end{equation}

 The following proposition holds
\begin{proposition}\label{Stokes operator propos}
Suppose that $\textbf{g}=(g',g_n) \in K^{l,\rho}_{\beta,T}$ with $\textbf{g}(t=0)=\textbf{0}$ and $g_n= |\xi'| \int \limits_0^{+\infty} dy' f(\xi',y',t) k(\xi',y')$ with $ |\xi'| \int \limits_0^{+\infty} d y' |k(\varepsilon',y')| \leq 1$ and $f \in L^{l,\rho,\theta}_{\beta,T}$. Then $\mathcal{S} \textbf{g} \in L^{l,\rho,\theta}_{\beta,T}$ and \begin{equation}
|\mathcal{S} \textbf{g}|_{l,\rho,\theta,\beta,T} \leq c ( |g'|_{l,\rho,\theta,\beta,T}+ |f|_{l,\rho,\theta,\beta,T} ).
\end{equation}
\end{proposition}
The proof can be found  in \cite{SC1998B} where the above results appears as   Proposition 3.4.

In appendix \ref{properties_of_PinftyE2}, we shall show that, when $\textbf{v} \in S^{l,\rho,1,\theta}_{\beta,T}$ then 
$\bar{P}^\infty \tilde{E}_2 \textbf{v} \in L^{l,\rho,\theta}_{\beta,T}$ and that, furthermore, 
$\gamma \bar{P}^\infty \tilde{E}_2$ satisfies the hypothesis of Proposition \ref{Stokes operator propos}; 
all this leads to the following result.
\begin{proposition}\label{NstarvinL}
Suppose that $\textbf{u} \in S^{l,\rho,1,\theta}_{\beta,T}$, then $\mathcal{N}^* \textbf{u} \in L^{l,\rho,\theta}_{\beta,T}$ and \begin{equation}
|\mathcal{N}^* \textbf{u}|_{l,\rho,\theta,\beta,T} \leq C |\textbf{u}|_{l,\rho,1,\theta,\beta,T}.
\end{equation}
\end{proposition}

We now give  some bounds  in a time integrated form that will be useful for the application of the ACK theorem. 
\begin{proposition}\label{Prop2}
Assume that $\textbf{u} \in S^{l,\rho,1,\theta}_{\beta,T}$. Then we have that, for $\rho' < \rho - \beta t$, $\theta' < \theta - \beta t$, \begin{equation}
|\mathcal{N}^* \textbf{u}|_{l,\rho',0,\theta'} \leq C \int \limits_0^t ds |\textbf{u}(\cdot,\cdot,s)|_{l,\rho',0,\theta'},
\end{equation}
\begin{equation}
|\partial_Y \mathcal{N}^* \textbf{u}|_{l-1,\rho',0,\theta'} \leq C \int \limits_0^t ds  \frac{|\textbf{u}(\cdot,\cdot,s)|_{l,\rho',1,\theta'}}{(t-s)^{1/4}}.
\end{equation}
%If $\gamma \textbf{u}=\textbf{0}$, then \begin{equation}
%|\mathcal{N}^* \textbf{u}|_{l,\rho',1,\theta'} \leq C \int \limits_0^t |\textbf{u}(\cdot,\cdot,s)|_{l,\rho',1,\theta'} + \varepsilon^{1/2} \frac{|\textbf{u}(\cdot,\cdot,s)|_{l-1/2,\rho',0,\theta'}}{(t-s)^{1/4}} ds.
%\end{equation}
\end{proposition}
Contrarily to the compatible case, the mild singularity in time cannot be completely eliminated, even when $\gamma \textbf{u} = \textbf{0}$: this is due to the fact that the functional setting $S^{l,\rho,1,\theta}_{\beta,T}$ is more singular than the one used in the compatible case. In particular, this setting does not allow a time derivative. The presence of this singularity implies that, in order to estimate a derivative of order $l+1$, like $\partial_Y \partial_x^l$, we cannot use the Cauchy estimates, otherwise we would have, at the denominator, $(t-s)^{1/4} (\rho(s)-\rho')$, which is not allowed by the ACK theorem. Therefore, we need to estimate the derivatives of order $l+1$ in a better way. A similar problem would appear if, in order to deal with the bilinear terms appearing in the error equation, one tries to use the algebra properties of $S^{l,\rho,1,\theta}_{\beta,T}$; therefore, we need some estimates in an $L^\infty_Y L^2_x$-like setting, so that the algebra properties are used only in the tangential variable. \\
The next proposition shows how to treat derivatives of order $l+1$
\begin{proposition}\label{Prop3}
For any $\textbf{u} \in S^{l,\rho,1,\theta}_{\beta,T}$, we have for $\rho' < \rho - \beta t$, $\theta' < \theta - \beta t$ 
\begin{equation}
|\mathcal{N}^* \textbf{u}|_{l+1,\rho',0,\theta'} \leq C \int \limits_0^t \frac{|\textbf{u}|_{l,\rho(s),0,\theta'}}{\rho(s)-\rho'} ds,
\end{equation} 
\begin{equation}
|\mathcal{N}^* \textbf{u}|_{l+1,\rho',0,\theta'} \leq \frac{C}{\varepsilon} \int \limits_0^t \frac{|\textbf{u}|_{l,\rho',0,\theta'}}{(t-s)^{1/2}} ds,
\end{equation}
\begin{equation}
|\partial_Y \mathcal{N}^* \textbf{u}|_{l,\rho',0,\theta'} 
\leq C  \int \limits_0^t \left[ \frac{|\textbf{u}|_{l,\rho',1,\theta'}}{(t-s)^{1/2}} + \frac{|\textbf{u}|_{l,\rho(s),1,\theta'}}{\rho(s)-\rho'} \right] ds,
\end{equation}
\begin{equation}
|\partial_{YY} \mathcal{N}^* \textbf{u}|_{l-1,\rho',0,\theta'} 
\leq C  \int \limits_0^t \left[  \frac{|\textbf{u}|_{l,\rho',1,\theta'}}{(t-s)^{1/2}} + \varepsilon^{1/2} \frac{|\textbf{u}|_{l-1/2,\rho',0,\theta'}}{(t-s)^{3/4}} \right]ds.
\end{equation}
Furthermore, if $\gamma \textbf{u}= \textbf{0}$, then \begin{equation}
|\partial_Y \mathcal{N}^* \textbf{u}|_{l,\rho',0,\theta'} \leq C  \int \limits_0^t \frac{|\textbf{u}|_{l,\rho',1,\theta'}}{(t-s)^{1/2}}  ds.
\end{equation}
\end{proposition}
The case $\gamma \textbf{u} = \textbf{0}$ is used to verify the quasi-contractiveness hypothesis of the ACK theorem: in this case, if we don't use Cauchy estimates, the only derivative of order $l+1$ which can be an $O\left(1/\varepsilon \right)$ is the purely tangential one, $\partial_x^{l+1} \mathcal{N}^* \textbf{u}$.

The next proposition provides some estimates in the $L^{\infty}_Y$ norm. 
\begin{proposition}\label{Prop4}
For any $\textbf{u} \in S^{l,\rho,1,\theta}_{\beta,T}$ we have that, $\forall \rho'<\rho-\beta t$, $\forall t \in [0,T]$ \begin{equation}
\sup \limits_{Y} |\mathcal{N}^* \textbf{u}|_{l,\rho'} \leq C \int \limits_0^t \frac{|\textbf{u}|_{l,\rho',0,\theta'}}{(t-s)^{1/4}}, 
\end{equation}
\begin{equation}
\sup \limits_{Y} |\partial_Y \mathcal{N}^* \textbf{u}|_{l-1,\rho'} \leq C \int \limits_0^t \frac{|\textbf{u}|_{l,\rho',0,\theta'}}{(t-s)^{3/4}}.
\end{equation}

%\begin{equation}
%\sup \limits_{Y} |\mathcal{N}^* \textbf{u}|_{l+1,\rho'} \leq \frac{C}{\varepsilon} \int \limits_0^t \frac{|\textbf{u}|_{l,\rho',0,\theta'}}{(t-s)^{3/4}} ,
%\end{equation}
%\begin{equation}
%\sup \limits_{Y} | \partial_Y \mathcal{N}^* \textbf{u}|_{l,\rho'} \leq C \int \limits_0^t \frac{|\textbf{u}|_{l,\rho',0,\theta'}}{(t-s)^{3/4}} ,
%\end{equation}
%\begin{equation}
%\sup \limits_{Y} | \partial_{YY} \mathcal{N}^* \textbf{u}|_{l-1,\rho'} \leq C |\textbf{u}|_{l,\rho,1,\theta,\beta,T}.
%\end{equation}

\end{proposition}

The proof of Propositions \ref{Prop2}, \ref{Prop3}, and \ref{Prop4}   can be found in Appendix \ref{proof_prop_7punto8910}. 

\subsection{The error equation}\label{sub::error_star}

Given the expression \eqref{edecompositione} for $\textbf{e}$, the system  \eqref{corrector_e} 
is now an equation for $\textbf{e}^*$ that can be cast in the following form:

\begin{equation}\label{estar}
\textbf{e}^*= \textbf{F}(\textbf{e}^*,t),
\end{equation}
where 
\begin{equation} \label{def::F}
\begin{split}
\textbf{F}(\textbf{e}^*,t)=\textbf{k} - \{[ \textbf{u}^{NS0} + \varepsilon (\textbf{u}^E_{(1)}+ \bar{\textbf{u}}^P_{(1)}+ \boldsymbol{\sigma} )] \cdot \nabla \mathcal{N}^* \textbf{e}^* + \mathcal{N}^* \textbf{e}^* \cdot \nabla [ \textbf{u}^{NS0} + \\ + \varepsilon (\textbf{u}^E_{(1)}+ \bar{\textbf{u}}^P_{(1)}+ \boldsymbol{\sigma})] + \varepsilon \mathcal{N}^*\textbf{e}^* \cdot \nabla \mathcal{N}^* \textbf{e}^* \}.
\end{split}
\end{equation}
The forcing term $\textbf{k}$ is given by 
\begin{equation}\label{kexpressionridotta}
\begin{split}
\textbf{k}= \boldsymbol{\Xi}_r+ \boldsymbol{\Psi}.
\end{split}
\end{equation}
In the above expression,  $\boldsymbol{\Xi}_r$ is the regular part of the source term appearing in the equation of the remainder $\textbf{e}$, 
see the system \eqref{corrector_e}. 
The expression for $\boldsymbol{\Xi}_r$ is reported in \eqref{Xi_r_expression}. 
The term   $\boldsymbol{\Psi}$ derives from the introduction of $\textbf{h}$ and $\boldsymbol{\sigma}$ in the decomposition \eqref{edecompositione};
the explicit expression of $\boldsymbol{\Psi}$ is given in \eqref{Psi_expression}. 
The following proposition states that the forcing term is bounded. 
\begin{proposition}\label{propoforcingterm}
There exists a constant $R_0$, independent of $\varepsilon$, such that \begin{equation}
|\textbf{k}|_{l,\bar{\rho},1,\bar{\theta},\bar{\beta},\bar{T}} \leq R_0
\end{equation}
\end{proposition} 
The key of the proof is to rearrange $\textbf{k}$ in a way such that terms that are $O(\varepsilon^{-1})$ multiply terms that are $O(\varepsilon)$, and terms that have singular derivative with respect to $Y$ multiply terms that go to zero with $Y$ or $t$. 
The details are reported  in the appendix \ref{proof_bounded_source}.

\subsubsection{Quasi contractiveness hypothesis}

In this subsection we shall prove that right hand side of equation \eqref{estar}, explicitly defined in \eqref{def::F}, 
satisfies the quasi-contractiveness hypothesis of the ACK Theorem. 
It is useful to distinguish between the nonlinear and the linear terms present in \eqref{def::F}. 

As for the nonlinear term, a special attention is needed for $ (\mathcal{N}^*(\textbf{e}^{*1} - \textbf{e}^{*2}))_n $ $ \partial_Y (\mathcal{N}^* \textbf{e}^{*2})'$: indeed, since $\gamma\textbf{e}^{*2} \neq \textbf{0}$, the $| \, |_{l,\rho',0,\theta'}$ norm of $ \partial_Y (\mathcal{N}^* \textbf{e}^{*2})'$ cannot be bounded by a constant, due to a term which behaves like $\partial_Y \tilde{E}_1 \gamma \bar{P}^{\infty} \tilde{E}_2 \textbf{e}^{2*}$. However, for the properties of $\tilde{E}_1$, $Y \partial_Y \tilde{E}_1 \gamma \bar{P}^{\infty} \tilde{E}_2 \textbf{e}^{2*}$ behaves like $\tilde{E}_1 \gamma \bar{P}^{\infty} \tilde{E}_2 \textbf{e}^{2*}$, while $(\mathcal{N}^*(\textbf{e}^{*1} - \textbf{e}^{*2}))_n /Y$ can be estimated in terms of $\partial_Y (\mathcal{N}^*(\textbf{e}^{*1} - \textbf{e}^{*2}))_n$. Therefore, multiplying and dividing by $Y$, we essentially move the normal derivative on the term with a more favourable estimate. \\ The other nonlinear terms are easily estimated using Propositions \ref{Prop2}, \ref{Prop3} and \ref{Prop4},  so we obtain the following result.

\begin{proposition}
	Suppose $\textbf{e}^{*i}\in S^{l.\rho',1,\theta'}$ for $i=1,2$. Then the following estimate holds:
\begin{equation} \label{primaquasicontr}
\left|\varepsilon \mathcal{N}^* \textbf{e}^{*1} \cdot \nabla \mathcal{N}^* \textbf{e}^{*1} - \varepsilon \mathcal{N}^* \textbf{e}^{*2} \cdot \nabla \mathcal{N}^* \textbf{e}^{*2} \right|_{l,\rho',1,\theta'} \leq C \int \limits_0^t \frac{|\textbf{e}^{*1} - \textbf{e}^{*2}|_{l,\rho',1,\theta'} }{(t-s)^{3/4}} ds.
\end{equation}
\end{proposition}
%\begin{equation}
%\varepsilon \mathcal{N}^* \textbf{e}^* \cdot \nabla \mathcal{N}^* \textbf{e}*= \varepsilon (\mathcal{N} \textbf{e}^*)_1 \partial_x \mathcal{N}^* \textbf{e}^* + (\mathcal{N} \textbf{e}^*)_2 \partial_Y \mathcal{N}^* \textbf{e}^*,
%\end{equation}  
%\begin{equation}
%\begin{split}
%\varepsilon \mathcal{N}^* \textbf{e}^{*1} \cdot \nabla \mathcal{N}^* \textbf{e}^{*1} - \varepsilon \mathcal{N}^* \textbf{e}^{*2} \cdot \nabla \mathcal{N}^* \textbf{e}^{*2}=  \\ \varepsilon \mathcal{N}^* (\textbf{e}^{*1} - \textbf{e}^{*2})  \cdot \nabla \mathcal{N}^* \textbf{e}^{*1} + \varepsilon \mathcal{N}^*  \textbf{e}^{*2} \cdot \nabla \mathcal{N}^* (\textbf{e}^{*1} - \textbf{e}^{*2}),
%\end{split}
%\end{equation} 
%so using 

For the linear terms, since $\gamma \textbf{e}^{*1}= \gamma \textbf{e}^{*2}$, we can apply the results of proposition \ref{Prop3} regarding the case with zero trace. The linear terms are problematic essentially for the presence of $u^S$ and its $Y$-derivatives, as well for the presence of the $O(\varepsilon^{-1})$ term 
$\partial_y \tilde{u}^P$. 
In Appendix \ref{sub::est_linear_F} we show how it is possible to estimate these terms. 

Proposition \ref{propoforcingterm}, the estimate of the nonlinear term \eqref{primaquasicontr}, and the estimates on the linear terms given in 
Appendix \ref{sub::est_linear_F}
 allow us to use the ACK theorem in the functional setting $S^{l,\bar{\rho},1,\bar{\theta}}_{\bar{\beta},T}$.
 This leads to the following result:
\begin{proposition}
Assume that $\textbf{u}_0 \in H^{l,\rho,\theta}$, with $\gamma_n \textbf{u}_0=0$; then  $ \textbf{e}^* \in S^{l,\bar{\rho},1,\bar{\theta}}_{\bar{\beta},T}$ for some $T>0$, $\bar{\rho}<\rho$, $\bar{\theta}<\theta$, $\bar{\beta}>\beta$, with all those parameters independent of the viscosity.
\end{proposition}
By proposition \ref{NstarvinL}, we therefore have $\mathcal{N}^* \textbf{e}^* \in L^{l,\bar{\rho},\bar{\theta}}_{\bar{\beta},T} $; combining with propositions \ref{hinL} and \ref{sigmainL}, we have that the overall error $\textbf{e}$ is in $L^{l,\bar{\rho},\bar{\theta}}_{\bar{\beta},T} $, and this concludes the proof of theorem \ref{mainth}.
\subsection{More general initial conditions}
The estimates performed for the forcing term in the equation of $\textbf{e}^*$ heavily relied on the assumption that the initial condition are purely eulerian, which means that $\textbf{e}(t=0)= \textbf{u}^E_{(1)}(t=0)= \bar{\textbf{u}}^P_{(1)}(t=0)= \bar{\textbf{u}}^P(t=0)=\textbf{0}$. The zero viscosity limit holds also for more general initial conditions like \begin{equation}
\textbf{u}_0= \textbf{u}^E_0 + \bar{\textbf{u}}^P_0 + \varepsilon \left(\textbf{u}^E_{(1),0}+ \bar{\textbf{u}}^P_{(1),0}+ \textbf{e}_{0} \right)
\end{equation}
with \begin{equation}
\gamma v^E_{(1),0}= - \gamma \bar{v}^P_{0}, \, \,  \gamma u^P_{(1),0}= -\gamma u^E_{(1),0}, \, \, \gamma \textbf{e}_0= (0,- \gamma \bar{v}^P_{(1),0}); 
\end{equation} this kind of initial data allows a zero order incompatibility with no-slip boundary condition, since we are not assuming that $\gamma \tilde{u}^P_0= - \gamma u^E_0$. The most challenging variation needed is to prove that the forcing term is still in $S^{l,\bar{\rho},1,\bar{\theta}}_{\bar{\beta},\bar{T}}$ and $O(1)$: we need to regroup some terms carefully in order to show the derivative $\partial_Y \textbf{k}$ has the desired regularity. This is done in appendix \ref{kgeneralinitial} \\

\section{Conclusions}\label{conclusionifinali}
In this paper, we have proved the short-time existence of analytic solutions of the Navier-Stokes equations in the half-space when 
the initial data are not compatible with the boundary data: the existence time of the solution we constructed is independent of the viscosity. 
When the viscosity goes to zero, the Navier-Stokes solution approaches the Euler solution away from the boundary, 
and the Prandtl solution inside a boundary layer whose thickness is proportional to the square root of the viscosity. 

Our method of constructing the Navier-Stokes solution is based on a matched asymptotic expansion: the solution is the superposition of an outer (Euler)
solution and an inner (Prandtl) solution, plus correction terms that are $O(\sqrt{\nu})$, see \eqref{informalmain} or \eqref{asymptotic_formal}. 
The inner Prandtl solution decays to zero exponentially outside the boundary layer, which gives the convergence of the Navier-Stokes solutions 
to the Euler solution in the $L^2$ norm. 

The incompatibility between the initial data and the no-slip condition causes the presence of singular terms which are $O(1)$ near the boundary and $O(\varepsilon)$ away from the boundary: the fact that the singularity is propagated in all the half space immediately through the term $ \varepsilon \sqrt{t} \boldsymbol{w}_b^S$ is coherent with the parabolic nature of the Navier-Stokes equations, which leads to an infinite speed of propagation of the singularity. Of course, when the viscosity goes to zero, the parabolic nature of the equations is lost, and indeed we have that $ \varepsilon \sqrt{t} \boldsymbol{w}_b^S$ goes to zero, leaving a singularity confined at the boundary.

We have imposed analytic initial data for two reasons: 
the first is that, in order to prove well posedness of the Prandtl equations without a monotonicity assumption, the use of data analytic in at least the tangential variable is "almost necessary"; 
the second reason is that we managed to obtain estimates independent of the viscosity by exploiting the restriction of the strip of analyticity in the tangential variable. 
\appendix

\section{The source terms} \label{appendixsource}

\subsection{The source term for the error equation \cref{corrector_e}}
The source term $\boldsymbol{\Xi}$ has the following expression:
\begin{equation} \label{Xi_soource}
	\begin{split}
\boldsymbol{\Xi} =\textbf{f}+\left(g\partial_y \tilde{u}^P,0\right) +\\
-\left[   \textbf{u}^{NS}_{(0)}   \cdot\nabla \bar{\textbf{u}}^{P}_{(1)}  +    \bar{\textbf{u}}^{P}_{(1)}  \cdot\nabla \textbf{u}^{NS}_{(0)}  +
\left(  \bar{\textbf{u}}^{P} +\varepsilon  \textbf{u}^{NS}_{(1)}  \right)\cdot \nabla  \textbf{u}^{E}_{(1)}+  
\textbf{u}^{E}_{(1)}  \cdot \nabla \left(  \bar{\textbf{u}}^{P} + \varepsilon  \bar{\textbf{u}}^{P}_{(1)}\right)   + \right. \\
\left. \varepsilon \bar{\textbf{u}}^{P}_{(1)}  \cdot \nabla \bar{\textbf{u}}^{P}_{(1)} \right] \\
+\varepsilon^2 \left[   \Delta {\textbf{u}}^{E}_{(1)} +\left(\partial_{xx} \bar{{u}}^{P}_{(1)},0\right) \right] -\left(0, \left(\partial_t-\varepsilon^2\Delta\right)
\bar{v}^P_{(1)}\right)
\end{split}
\end{equation}
where $\textbf{f}$ is given by: 
\begin{equation}\label{fregsing}
	\textbf{f}= \textbf{f}_r+ (\varepsilon \partial_x^2 u^S,0),
\end{equation}
\begin{equation}\label{freg}
	\begin{split}
		f^1_r= - \varepsilon^{-1} \{ \tilde{u}^P ( \partial_x u^E - \partial_x u^E\big|_{y=0}) + \partial_x \tilde{u}^P (u^E-u^E\big|_{y=0})+ \\ \partial_y \tilde{u}^P (v^E + y \partial_x u^E\big|_{y=0}) \}
		-\bar{v}^P \partial_y u^E + \varepsilon \Delta u^E + \varepsilon \partial_x^2 \tilde{u}^R, \\
		f^2_r=- [ \partial_t \bar{v}^R + u^{NS0} \partial_x \bar{v}^P + v^{NS0} \partial_y \bar{v}^P + \bar{v}^P \partial_y v^E ] - \varepsilon^{-1} \tilde{u}^P \partial_x v^E + \\ + \varepsilon \Delta v^E + \varepsilon^2 \Delta v^R + \varepsilon^2 \partial_x^2 \tilde{v}^S.
	\end{split}
\end{equation}
Writing down the expression of $f^2$, we used the fact that $\partial_t\bar{v}^S - \partial_Y^2 \bar{v}^S=0$ to get rid of those singular terms. We can write $g=g^R+g^S$, with \begin{equation}
	g^S= - \partial_x u_0(x,y=0) \int \limits_0^{+ \infty} d Y \int \limits_Y^{+ \infty} \frac{e^{\frac{- Y'^2}{4 t}}}{\sqrt{t \pi}} d Y' = C_1 \sqrt{t} \partial_x u_0(x,y=0) . 
\end{equation}
In equation \eqref{fregsing}, we isolated the most singular term of $\textbf{f}$; in a similar way, we write 
\begin{equation}\label{Xi_r_expression}
\boldsymbol{\Xi}_r=\boldsymbol{\Xi}-\varepsilon (\partial_x^2 u^S,0),
\end{equation}
where $\boldsymbol{\Xi}$ is given in \eqref{Xi_soource}. Notice that $\textbf{f}+\left(g\partial_y \tilde{u}^P,0\right)$ is the forcing term in the equation of the remainder $\textbf{R}$ in the decomposition \begin{equation}
\textbf{u}^{NS}=\textbf{u}^E + \bar{\textbf{u}}^P+ \varepsilon \textbf{R}
\end{equation}
The term $\textbf{f}$ is $O(1)$ with respect to $\varepsilon$, while $g \partial_y \tilde{u}^P$ is $O( \varepsilon^{-1})$: so an additional decomposition of the remainder $\textbf{R}$ is needed to prove the validity of the inviscid limit.  

\subsection{Estimates for the forcing term of $\textbf{e}^*$}  \label{proof_bounded_source}
The source term $\textbf{k}$ in  \eqref{kexpressionridotta} is decomposed as $\textbf{k}= \boldsymbol{\Xi}_r+\boldsymbol{\Psi}$ where 

 \begin{equation}\label{Psi_expression}
\begin{split}
\boldsymbol{\Psi}= - \bigg((\textbf{u}^{NS0}+\varepsilon(\textbf{u}^E_{(1)}+\bar{\textbf{u}}^P_{(1)})) \cdot \nabla( \textbf{h}+ \boldsymbol{\sigma}) + (\textbf{h} + \boldsymbol{\sigma}) \cdot \nabla (\textbf{u}^{NS0}+\varepsilon(\textbf{u}^E_{(1)}+\bar{\textbf{u}}^P_{(1)})) \\ + (\textbf{h} + \boldsymbol{\sigma}) \cdot \nabla(\textbf{h} + \boldsymbol{\sigma}) \bigg) + \varepsilon^2 \partial_{xx} (\textbf{h} + \boldsymbol{\sigma})-\varepsilon(0, (\partial_t- \partial_{YY}) h_n)= \\
 - \{[(\textbf{u}^{NS0} + \varepsilon (\textbf{u}^E_{(1)}+ \bar{\textbf{u}}^P_{(1)}+\textbf{h}+ \boldsymbol{\sigma})) \cdot \nabla (\bar{ \textbf{u}}^P_{(1)}+\textbf{h})+ (\textbf{u}^E_{(1)} \cdot \nabla \bar{\textbf{u}}^P- (g \partial_y \tilde{u}^P,0))+ \\ ( \bar{\textbf{u}}^P_{(1)}+\textbf{h}+ \boldsymbol{\sigma}) \cdot 
  \nabla \bar{\textbf{u}}^P ] + [( \textbf{u}^E_{(1)}+ \bar{\textbf{u}}^P_{(1)}+\textbf{h}) \cdot \nabla \textbf{u}^E + (\bar{u}^P+ \varepsilon(\textbf{u}^E_{(1)}+ \bar{\textbf{u}}^P_{(1)}+\textbf{h}+ \boldsymbol{\sigma})) \cdot \\ 
  \nabla \textbf{u}^E_{(1)} +( \textbf{u}^{NS0} + \varepsilon (\textbf{u}^E_{(1)}+ \bar{\textbf{u}}^P_{(1)}+\textbf{h}+ \boldsymbol{\sigma} )) \cdot  \nabla \boldsymbol{\sigma}]\}+ 
  \textbf{f}_r - (0,(\partial_t - \varepsilon^2 \Delta)\varepsilon (\bar{v}^P_{(1)}+h_n))
    \\ + \varepsilon^2[ \Delta \textbf{u}^E_{(1)} + \partial_{xx} (u^P_{(1)}+h',0) + \partial_{xx} \boldsymbol{\sigma}].
\end{split}
\end{equation}
We want to prove that the forcing term $\textbf{k}$  is in $S^{l,\bar{\rho},1,\bar{\theta}}_{\bar{\beta},\bar{T}}$ and 
that it is $O(1)$ with respect to $\varepsilon$. In our analysis we shall focus on the most challenging terms.  
We begin with:
\begin{equation}
[\textbf{u}^{NS0} + \varepsilon (\textbf{u}^E_{(1)}+ \bar{\textbf{u}}^P_{(1)}+\textbf{h}+ \boldsymbol{\sigma})) ] \cdot \nabla (\bar{ \textbf{u}}^P_{(1)}+\textbf{h}) .
\end{equation} 
This term is $O(1)$ because $\partial_y (\bar{ \textbf{u}}^P_{(1)}+\textbf{h})$, which has an $O(1/\varepsilon)$ $L^2_Y$ norm, multiplies a term with an $O(\varepsilon)$ $L^\infty_Y$ norm. When we take a partial derivative with respect to $Y$, some singular terms appear, $\partial_Y \tilde{u}^S \partial_x (\bar{ \textbf{u}}^P_{(1)}+\textbf{h})$ and $ [\textbf{u}^{NS0}/\varepsilon + (\textbf{u}^E_{(1)}+\bar{ \textbf{u}}^P_{(1)}+\textbf{h}+\boldsymbol{\sigma}) ]_2 \partial_{YY} u^S_1$: for the first term, since $v^E= \varepsilon Y (\partial_y v^E)(\lambda(y) y)$ for some $\lambda(y) \in [0,1]$, we have that $\bar{ \textbf{u}}^S_{(1)}$ goes to zero as $\sqrt{t}$ in the $L^\infty_Y$ norm, $\bar{ \textbf{u}}^R_{(1)}$ and $\textbf{h}$ go to zero linearly with $t$, and this is enough to balance the singularity of $\partial_Y \tilde{u}^S$; for the second term,  $v^E/\varepsilon+ \bar{v}^S+(\textbf{u}^E_{(1)}+ \bar{\textbf{u}}^P_{(1)}+\textbf{h}+ \boldsymbol{\sigma}))_2$ goes to zero linearly with $Y$, and this allows to balance the singularity of $\partial_{YY} u^S_{(1)}$. 
So all the singular terms in $[\textbf{u}^{NS0} + \varepsilon (\boldsymbol{w}+ \boldsymbol{w} + \boldsymbol{\sigma} ) ] \cdot \nabla (\bar{ \textbf{u}}^P_{(1)}+\textbf{h})$ are balanced, and its norm is controlled through the usual arguments involving algebra properties. 

We can now pass to estimate the terms: 
\begin{equation}
\textbf{u}^E_{(1)} \cdot \nabla \bar{\textbf{u}}^P- (g \partial_y \tilde{u}^P,0)
\end{equation}
The most problematic term is $(v^E_{(1)}-g) \partial_y \tilde{u}^P$, which is $O(1)$ since $v^E_{(1)}-g$ goes to zero linearly in $y=\varepsilon Y$; since $\partial_y u^S_{(1)}$ is multiplied by $\varepsilon Y$ times a bounded function, the product has finite $|\cdot|_{l,\rho',0,\theta'}$ norm. When we take the partial derivative with respect to $Y$, $\partial_Y v^E_{(1)}= \varepsilon \partial_y v^E_{(1)}$, so $\partial_Y v^E_{(1)} \partial_y \tilde{u}^P$ is still $O(1)$, and since $\partial_Y v^E_{(1)}$ has a regular part which goes to zero linearly with $t$ and a singular part which goes to zero with $t$ like $\sqrt{t}$, this product has finite $|\cdot|_{l,\rho',0,\theta'}$ norm. Finally, $(v^E_{(1)}-g) \partial_{yY} \tilde{u}^P= (\partial_Y v^E_{(1)})(\lambda y)Y \partial_{YY} \tilde{u}^P$, with the regular part of $\partial_Y v^E_{(1)}$ going to zero linearly with $t$, while the singular part goes to zero like $\sqrt{t}$. 

We  now  estimate the terms: 
\begin{equation}
(\bar{\textbf{u}}^P_{(1)} + \textbf{h}+ \boldsymbol{\sigma}) \cdot \nabla \bar{\textbf{u}}^P.
\end{equation}
Both $(u^P_{(1)} + h' + \sigma_1) \partial_Y \partial_x \tilde{u}^P$ and $ \left( \varepsilon\partial_Y\bar{v}^P_{(1)} + \varepsilon \partial_Y h_n + \partial_Y \sigma_2/\right) \partial_y \bar{\textbf{u}}^P =- (\partial_x u^P_{(1)} + \partial_x h' + \partial_x \sigma_1) \partial_Y \bar{\textbf{u}}^P$ are products of terms whose $L^\infty$ norm is $O(1)$ and goes to zero with t at least as $\sqrt{t}$ and a function whose $L^2_Y$ behaves like $t^{-1/4}$; in $( \bar{v}^P + h_n + \sigma_2/\varepsilon) \partial_{YY}\bar{\textbf{u}}^P$, the first factor has an $O(1)$ $L^\infty_Y$ norm which goes to zero linearly with t, and this balances the singularity of the second factor.

The term
\begin{equation}
\varepsilon^{-1} \tilde{u}^P ( \partial_x u^E - \partial_x \gamma u^E)
\end{equation}
 is $O(1)$ because $\partial_x u^E - \partial_x \gamma u^E$ goes to zero linearly with $y= \varepsilon Y$, and for the same reason $\varepsilon^{-1} \partial_Y \tilde{u}^P ( \partial_x u^E - \partial_x \gamma u^E)$ is $O(1)$ and with finite $|\cdot|_{l,\rho',0,\theta'}$ norm, while $\varepsilon^{-1} \tilde{u}^P \partial_Y \partial_x u^E=\tilde{u}^P \partial_y \partial_x u^E$ is of course $O(1)$ and bounded. The term $\varepsilon^{-1} \partial_x\tilde{u}^P ( u^E - \gamma u^E)$ can be treated in a similar way.

Finally, the term
\begin{equation}
\frac{1}{\varepsilon}\partial_y \tilde{u}^P (v^E + y \partial_x \gamma u^E)=\frac{1}{\varepsilon^2} \partial_Y \tilde{u}^P (v^E + y \partial_x \gamma u^E)
\end{equation}
is $O(1)$ because $(v^E + y \partial_x \gamma u^E)$ goes to zero quadratically in $y^2=\varepsilon^2 Y^2$ due to the incompressibility condition, and for the same reason 
 $ \frac{1}{ \varepsilon^2} \partial_{YY} \tilde{u}^P (v^E+ y \partial_x \gamma u^E)$ is $O(1)$ and with finite
  $| \cdot |_{l,\rho',0,\theta'}$ norm, while
  $\frac{1}{\varepsilon} \partial_{Y} \tilde{u}^P (\partial_y v^E+  \partial_x \gamma u^E)$ has $\partial_y v^E+  \partial_x \gamma u^E$ which goes to zero linearly with $\varepsilon Y$
  
All the other terms in $\textbf{k}$ are easier to estimate.

\subsection{Forcing term \textbf{k} with general initial conditions}\label{kgeneralinitial} 
When the initial conditions are not purely eurelian, we arrange the terms so that their sum is zero at the boundary. 
This is done in a different way for the tangential and the normal components.

The sum of $\partial_Y u^S \partial_x w_1^{R}$, deriving from the derivative of $u^S \partial_x w_1^{R}$, and $\partial_Y u^S \partial_x u_{(1)}^{R}$, deriving from the derivative of $u^S \partial_x u_{(1)}^{R}$, has finite $|\cdot|_{l-1,\rho',0,\theta'}$ norm because $\partial_x w_1^{R} + \partial_x u_{(1)}^{R}$ goes to zero linearly in $Y$ 
(as well in $t$, since $\gamma w_{1,0}^{RR}=-\gamma w_{1,0}^{RR}$). 

The sum of  $\partial_Y u^S \varepsilon \partial_x \bar{v}_{(1)}^{R}$, deriving from the derivative of $u^S \partial_x \bar{v}_{(1)}^{R}$, and $\partial_Y \tilde{u}^S \partial_x \sigma_2$, deriving from the derivative of $\tilde{u}^S \partial_x \sigma_2$, has finite $|\cdot|_{l-1,\rho',0,\theta'}$ norm because $\varepsilon \partial_x \bar{v}_{(1)}^{R}+ \partial_x \sigma_2$ goes to zero linearly in $Y$ (as well in $t$). 

The sum of $w_1^{R} \partial_x \partial_Y u^S$, deriving from the derivative of $w_1^{R} \partial_x  u^S$, and $ u_{(1)}^{R} \partial_x \partial_Y u^S$, deriving from the derivative of $w_1^{RR} \partial_x u^S$, has finite $|\cdot|_{l-1,\rho',0,\theta'}$ norm because  $ w_1^{R} +  u_{(1)}^{R}$ goes to zero linearly in $Y$ (and also in $t$). 

The sum of $\partial_y w_2^{R} \partial_Y u^S=- \partial_x w_1^{R} \partial_Y u^S$, 
deriving from the derivative of $(w_2^{R} - g^R) \partial_y u^S$, 
and $(\partial_Y \bar{v}_{(1)}^{R} + \varepsilon^{-1} \partial_Y \sigma_2) \partial_Y u^S= -( \partial_x u_{(1)}^{R} + \partial_x \sigma_1) \partial_Y u^S$, 
deriving from the derivative of $( \bar{v}_{(1)}^{R} + \varepsilon^{-1} \sigma_2) \partial_Y u^S$, 
has finite $|\cdot|_{l-1,\rho',0,\theta'}$ norm because 
$\partial_x w_1^{R}+ \partial_x u_{(1)}^{R}$ 
goes to zero linearly in $Y$ (and in $t$) and  
$\partial_x \sigma_1$ goes  to zero linearly with $t$. 

The sum of $-w_2^{R} \partial_{YY} \bar{v}^S$, deriving from the derivative of $-w_2^{RR} \partial_{Y} \bar{v}^S$, and $\varepsilon^{-1} v^{NS0}$ $ \partial_{YY} \bar{v}^S$, deriving from the derivative of $\varepsilon v^{NS0} \partial_{Y} \bar{v}^S$, has finite $|\cdot|_{l-1,\rho',0,\theta'}$ norm: indeed, $v^{NS0}= v^E + \varepsilon \bar{v}^R + \varepsilon \bar{v}^S$, where $v^E$ goes to zero linearly with $y=\varepsilon Y$, $\bar{v}^S$ goes to zero like $\sqrt{t}$ and $-w_2^{R}+ \bar{v}^R$ goes to zero linearly in $Y$ (and in $t$). 

For all the other terms there is no need to change the argument used in the case $\textbf{u}_0= \textbf{u}_0^E$.

\section{Proof of Proposition \ref{prop:estimate_h}}  \label{proof_on_h}

The estimate on $h'$ is a direct consequence of the expression \eqref{expression_h_prime}, and of the following estimates on $F$ and on its $Y$-derivatives. 
We begin with $F$. 

\begin{eqnarray} 
		\sup \limits_{Y \in \Sigma(\theta)} e^{\mu Y} \left|  F \right|= \sup \limits_{Y} \left| \frac{2}{\sqrt{\pi}} \int \limits_0^t ds \int \limits_0^{+\infty} e^{\mu Re(Y-Y')}(E_0^--E_0^+) e^{\mu Re(Y')}\int \limits_{\frac{Y'}{2\sqrt{t}}}^{+\infty} e^{-\sigma^2} d \sigma d Y' \right| \leq  \nonumber \\
\sup \limits_{Y} \left| \frac{2}{\sqrt{\pi}} \int \limits_0^t ds \int \limits_0^{+\infty} e^{\mu Re(Y-Y')}(E_0^--E_0^+) \right| \leq 
c t \leq C\nonumber
\end{eqnarray}
To pass from the first to the second line we have taken  into account the argument of subsection \ref{pathsubsec} about the complex gaussian 
and the fact that
$$	
\int \limits_{\frac{Y'}{2\sqrt{s}}}^{+\infty} e^{-\sigma^2} d \sigma \leq e^{-\frac{Y'^2}{8s}}\int \limits_{0}^{+\infty} e^{-\frac{\sigma^2}{2}} d \sigma,
$$
while, to get the last inequality, we have used that for $z \in \mathbb{R}$ 
$$	
e^{\mu z} e^{-c\frac{z^2}{t}} \leq e^{\mu z} e^{-c\frac{z^2}{2T}} e^{-c\frac{z^2}{2t}} \leq C(\mu,T) e^{-c\frac{z^2}{2t}}.
$$

A similar argument shows that 
$$
	\sup \limits_{Y \in \Sigma(\theta)} e^{\mu Y} \left| \partial_Y F \right | \leq c\sqrt{t}.
$$

To estimate  $\partial_{YY}F$, we first compute $\partial_{Y}$. Integrating by parts, we have 
\begin{equation}
	\partial_Y F= -  \int \limits_0^t \frac{e^{-\frac{Y^2}{4(t-s)}}}{\sqrt{ \pi (t-s)}}ds + \int \limits_0^t ds \int \limits_0^{+\infty} (E_0^-+E_0^+) \frac{e^{-\frac{Y'^2}{4s}}}{\sqrt{\pi s}} d Y'.
\end{equation}
We can now estimate $\partial_{YY}F$ as follows:
\begin{eqnarray}
e^{\mu Re(Y)}	|\partial_{YY} F|= 
\nonumber\\
e^{\mu Re(Y)}\left| \int \limits_0^t \frac{Y e^{-\frac{Y^2}{4(t-s)}}}{2\sqrt{\pi}(t-s)^{3/2}} ds + \int \limits_0^t ds \int \limits_0^{+\infty} \partial_Y(E_0^-+E_0^+) \frac{e^{-\frac{Y'^2}{4s}}}{\sqrt{\pi s}} d Y'\right|\leq 
\nonumber\\
e^{\mu Re(Y)}  \int \limits_{\frac{Y}{2\sqrt{t}}}^{+\infty} \frac{e^{-\sigma^2}}{\sqrt{\pi}} d \sigma + e^{\mu Re(Y-Y')} +  
\nonumber \\
 \int \limits_0^t ds \int \limits_0^{+\infty}  e^{\mu Re(Y-Y')}  \partial_Y(E_0^-+E_0^+)   e^{\mu Re(Y')}  \frac{e^{-\frac{Y'^2}{4s}}}{\sqrt{\pi s}} d Y'\leq
\nonumber \\
C+C \int \limits_0^t \frac{ds}{\sqrt{s(t-s)}} \leq C \nonumber 
\end{eqnarray}
In the first of the inequalities above, in the first integral we have simply used a standard change  of variable. 
To pass from the third line to the fourth, in the first integral we have used the exponential decay of the integral while, 
in the second integral, the fact that the gaussian in $Y'$ dominates the exponential $e^{\mu Re(Y')}$,  that 
the gaussian in $Y-Y'$ dominates the exponential $e^{\mu Re(Y-Y')}$, and a standard change of variable. 

The above estimate concludes the bound on $h'$. 
The estimate on the normal component $h_n$ is a direct consequence 
of the expression \eqref{expression_h_n}, and of the above estimates on $F$. 
One could in fact say more:

\begin{equation}
	\sup \limits_{Y \in \Sigma{(\theta)}} | \partial_{Y}^j h_n |_{l-3-j,\rho} \leq C \varepsilon t^{\frac{3-j}{2}} |u_{00}|_{l-j,\rho}.
\end{equation}

\section{Properties of $\bar{P}^\infty \tilde{E}_2$}\label{properties_of_PinftyE2}

It's easy to see that \begin{equation}\label{BoundednessPinfbar}
|\bar{P}^\infty\textbf{v}|_{l,\rho,j,\theta} \leq C |\textbf{v}|_{l,\rho,j,\theta}.
\end{equation}
Indeed, the convolution part of $\bar{P}^\infty\textbf{u}$ can be treated with Young's inequality for convolutions, with the $L^1$ norm of $\varepsilon|\xi'| e^{-|\xi'|\varepsilon Y} \chi_{Y\geq 0}$ equal to 1; furthermore,  taking derivatives with respect to Y, when the derivatives hit the extreme of integration the exponentials multiplying  $u_1$ and $u_2$ are bounded and the $\varepsilon |\xi'|$ outside the integral is balanced by the fact that we are taking fewer derivatives with respect to $x$; on the other side, when the derivatives hit the exponentials inside the integral, this simply causes the appearance of an $\varepsilon |\xi'|$. and 

The same argument shows that
 \begin{equation}\label{imprecise}
| \sup \limits_Y |\partial_Y^j \partial_x^i \bar{P}^{\infty} \textbf{v}| |_{0,\rho} \leq C \sum \limits_{h+k \leq i+j, h \leq i} | \sup \limits_Y |\partial_Y^h \partial_x^{k}\textbf{v}| |_{0,\rho}.
\end{equation}
Moreover, since only the non convolution part of $\bar{P}^{\infty'}$ actually sees all the partial derivatives with respect to $Y$, we can write 
\begin{equation} \label{precise}
 \sup \limits_Y |  \partial_Y^j \partial_x^i \bar{P}^{\infty} \textbf{v} |_{0,\rho} \leq \sup \limits_Y | \partial_Y^j \partial_x^i  v_1 |_{0,\rho} + C \sum \limits_{h+k \leq i+j, h < j} | \sup \limits_Y |\partial_Y^h \partial_x^{k}\textbf{v}| |_{0,\rho}.
\end{equation}
The above two estimates imply that, instead of $|\bar{P}^{\infty} \tilde{E}_2 \textbf{u}|_{l,\rho,2,\theta}$,  we can estimate $| \tilde{E}_2 \textbf{u}|_{l,\rho,2,\theta}$, and that, instead of $\sup \limits_Y | \partial_Y^i \partial_x^j \bar{P}^{\infty} \tilde{E}_2 \textbf{u}|_{0,\rho}$,
 we can estimate the right hand side of \eqref{precise}. 

It is also useful to notice that
\begin{equation}
\sup \limits_{Y \in \Gamma(\theta')} |f|_{0,\rho} \leq |\sup \limits_{Y \in \Gamma(\theta')} | f| |_{0,\rho} \label{inversion}
\end{equation}
that  can be proven  using Minkowski's integral inequality.

We now pass to the estimates involving the $\tilde{E}_2$ operator. 
Using Minkowski integral inequality to pass under integral sign, then Young inequality for convolutions in both $x$ and $Y$, we have that \begin{equation}\label{xL2normal}
|\partial_x^i \tilde{E}_2 f|_{0,\rho',\theta'}(t) \leq C \int \limits_0^t ds |\partial^i_x f|_{0,\rho',\theta'}(s) \leq C  |f|_{l,\rho,0,\theta,\beta,T} 
\end{equation}
$\forall i=0,...,l, \, \, \, \, \forall t, \, \, \, \forall \rho' \leq \rho-\beta t, \, \, \, \, \forall \theta' \leq \theta - \beta t$. If we want to take an extra partial derivative, we can either use Cauchy estimates or we can let the derivative hit the kernel, the cost being the appearance of an unbalanced $\frac{1}{\varepsilon\sqrt{t-s}}$ \begin{equation}\label{xL2extrarho}
|\tilde{E}_2 f|_{l+1,\rho',0,\theta'} \leq C \int \limits_0^t \frac{|f|_{l,\rho(s),0,\theta'}}{\rho(s)-\rho'} ds,
\end{equation} 
\begin{equation}\label{xL2extrat}
|\tilde{E}_2 f|_{l+1,\rho',0,\theta'} \leq \frac{C}{\varepsilon} \int \limits_0^t \frac{|f|_{l,\rho',0,\theta'}}{\sqrt{t-s}} ds.
\end{equation} 
When we use an $L^\infty$ norm for $Y \in \Gamma(\theta')$, we still use Young's inequality for convolution, but this time we use the $L^2$ norm (with respect to $Y$) on the kernel, so an additional $\frac{1}{(t-s)^{1/4}}$ appears: we have that \begin{equation}\label{xLinfnormal}
\begin{split}
| \sup \limits_{Y \in \Gamma(\theta')}|\partial_x^i \tilde{E}_2 f||_{0,\rho} \leq C \int \limits_0^t  \frac{| |\partial_x^if|_{L^2(\gamma(\theta'))}|_{0,\rho}}{(t-s)^{1/4}} ds= C \int \limits_0^t  \frac{| |\partial_x^i f|_{0,\rho}|_{L^2(\gamma(\theta'))}}{(t-s)^{1/4}} ds \\ \leq C |f|_{l,\rho,0,\theta,\beta,T}.
\end{split}
\end{equation}
The key point here is that the $L^2$ norm in the $Y$ variable and the Hardy norm in the $x$ variable are interchangeable, essentially because the Hardy norm is an $L^2$ norm on the boundary lines. We also obtain that \begin{equation}\label{xLinfextra}
| \sup \limits_{Y \in \Gamma(\theta')}| \partial_x^{l+1} \tilde{E}_2 f||_{l+1,\rho} \leq \frac{C}{\varepsilon} \int \limits_0^t  \frac{| |f|_{l,\rho}|_{L^2(\gamma(\theta'))}}{(t-s)^{3/4}} ds.
\end{equation}
The same argument works for the first partial derivative with respect to $Y$: while we make the derivatives with respect to $x$ act on $f$, we can leave the derivative with respect to $Y$ on the kernel, so we obtain \begin{equation}\label{YL2norm}
|\partial_x^i \partial_Y \tilde{E}_2 f|_{0,\rho',\theta'} \leq C \int \limits_0^t ds  \frac{|\partial_x^i f|_{0,\rho',\theta'}(s)}{\sqrt{t-s}}  \leq C |f|_{l,\rho,0,\theta,\beta,T}
\end{equation}
$\forall i=0,...,l-1, \, \, \, \, \forall t, \, \, \, \forall \rho' \leq \rho-\beta t, \, \, \, \, \forall \theta' \leq \theta - \beta t$. If $\gamma f=0$, then we can use integration by parts to move $\partial_Y=-\partial_{Y'}$ on $f$, so we obtain \begin{equation}\label{YL2normgamma}
|\partial_x^i \partial_Y \tilde{E}_2 f|_{0,\rho',\theta'} \leq C \int \limits_0^t |\partial_x^i \partial_Y f|_{0,\rho',\theta'}(s) ds.
\end{equation}
If we want to take $l$ derivatives with respect to $x$, we still have \begin{equation}\label{YL2extra}
|\partial_Y \tilde{E}_2 f|_{l,\rho',0,\theta'} \leq C \int \limits_0^t ds  \frac{|f|_{l,\rho',0,\theta'}(s)}{\sqrt{t-s}}.
\end{equation}
Using an $L^\infty$ norm on $Y$, we have that 
\begin{equation}\label{YLinfnormeextra}
|\sup \limits_{Y \in \Gamma(\theta')}| \partial_x^i \partial_Y \tilde{E}_2||_{0,\rho'} \leq C \int \limits_0^t \frac{|f|_{l,\rho',0,\theta'}}{(t-s)^{3/4}}ds.
\end{equation}
$i =0,...,l$
Before we take a second partial derivative with respect to $Y$, we first perform an integration by parts; if $\gamma f=0$, then \begin{equation}\label{YYL2normeextragamma}
|\partial_{YY} \tilde{E}_2 f|_{i,\rho',0,\theta'} \leq C \int_0^t ds  \frac{|\partial_Yf|_{i,\rho',0,\theta'}}{(t-s)^{1/2}}
\end{equation}
$i=0,...,l-1$, while if $\gamma f \neq 0$ we have an additional term given by \begin{equation}\label{YYL2normeextra}
\begin{split}
C \left| \left| \int \limits_0^t ds \frac{Y e^{-\frac{Y^2}{4(t-s)}}}{2 (t-s)^{3/2}} \gamma \partial_x^i \tilde{f} ds \right|_{L^2(\Gamma(\theta'))} \right|_{0,\rho'} \leq C \left| \int \limits_0^t \frac{|\partial_x^i \tilde{f}|_{W^{1,2}(\Gamma(\theta'))}}{(t-s)^{3/4}} ds \right|_{0,\rho'} \\ \leq  C \int \limits_0^t \frac{| f|_{l,\rho',0,\theta'}(s)}{(t-s)^{3/4}} ds ,
\end{split}
\end{equation} where \begin{equation}
\tilde{f}= \int \limits_{-\infty}^{+\infty} \frac{e^{-\frac{(x-x')^2}{4\varepsilon^2 (t-s)}}}{\sqrt{4 \pi (t-s)}} f(s,x) dx'.
\end{equation} 
In the $L^\infty$ norm, if $\gamma f=0$ then \begin{equation}\label{YYLinftynormeextragamma}
||\partial_x^i \partial_{YY} \tilde{E}_2 f|_{L^\infty(\Gamma(\theta'))} |_{0,\rho'} \leq C \int \limits_0^t \frac{|\partial_Y f|_{i,\rho',0,\theta'}}{(t-s)^{3/4}} ds ,
\end{equation}
$i=0,...,l-1$; otherwise, since we have an additional term given by \begin{equation}\label{boundary}
C  \int \limits_0^t \frac{Y e^{-\frac{Y^2}{4(t-s)}}}{2 (t-s)^{3/2}} \gamma \partial_x^i \tilde{f} ds 
\end{equation}
this term has finite $||\cdot|_{l,\rho}|_{L^\infty(\Sigma(\theta'))}$ norm, because \begin{equation}\label{YYLinftynormeextra}
\begin{split}
C \left| \left|  \int \limits_0^t \frac{Y e^{-\frac{Y^2}{4(t-s)}}}{2 (t-s)^{3/2}} \gamma \partial_x^i \tilde{f} ds \right|_{0,\rho'} \right|_{L^\infty_Y} \leq C  \left|  \int \limits_0^t  \frac{Y e^{-\frac{Y^2}{4(t-s)}}}{2 (t-s)^{3/2}} \gamma |\partial_x^i f|_{0,\rho'} ds  \right|_{L^\infty_Y} \\ \leq C \left|  \int \limits_0^t  \frac{Y e^{-\frac{Y^2}{4(t-s)}}}{2 (t-s)^{3/2}} | f|_{l,\rho',1,\theta'} ds  \right|_{L^\infty_Y} \leq C |f|_{l,\rho,1,\theta,\beta,T} \sup \limits_Y |\int \limits_{\frac{Y}{2\sqrt{t}}}^{+\infty} e^{-\sigma^2} d \sigma| \\ \leq C |f|_{l,\rho,1,\theta,\beta,T}.
\end{split}
\end{equation}
Notice that, in this case, we are not able to use the $||\cdot|_{L^\infty(\Sigma(\theta'))}|_{l,\rho}$ norm, because either we have to take the supremum with respect to $s$ before we take the $|\cdot|_{0,\rho'}$ norm (in order to take $f$ outside the time integral), or we have to bring the $L^{\infty}_Y$ norm inside the time integral, with $\sup \limits_Y |\frac{Y e^{-\frac{Y^2}{4(t-s)}}}{2 (t-s)^{3/2}}| = \frac{C}{t-s}$, which is not integrable in time.

\section{Proof of  \cref{NstarvinL}}
The results of the previous section, 
and the obvious equality $\partial_t \tilde{E}_2 \textbf{u}=\textbf{u} + \varepsilon^2 \partial_{xx} \tilde{E}_2 \textbf{u} + \partial_{YY} \tilde{E}_2\textbf{u}$, 
show that $\bar{P}^\infty \tilde{E}_2 \textbf{u} \in L^{l,\rho,\theta}_{\beta,T}$, with
\begin{equation}
|\bar{P}^\infty \tilde{E}_2 \textbf{u}|_{l,\rho,\theta,\beta,T} \leq C |\textbf{u}|_{l,\rho,1,\theta,\beta,T}.
\end{equation}
Now we want to prove that the trace of $ \bar{P}^{\infty} \tilde{E}_2 \textbf{u}$ is regular enough for the application of proposition 3.3 in \cite{SC1998B}. 
We have \begin{equation}\label{Pbarnormbord}
\gamma \bar{P}^\infty_n (\tilde{E}_2 \textbf{u})= \varepsilon |\xi'| \int \limits_0^{+ \infty} e^{- \varepsilon |\xi'| Y'} \frac{i \xi' }{|\xi'|} \tilde{E}_2 u_1 d Y',
\end{equation}
\begin{equation}\label{Pbartanbord}
\gamma \bar{P}^{\infty '} (\tilde{E}_2 \textbf{u})= \varepsilon |\xi'| \int \limits_0^{+ \infty} e^{- \varepsilon |\xi'| Y'} \frac{i \xi' }{|\xi'|} \tilde{E}_2 u_2 d Y'.
\end{equation}
We shall prove that, if $\textbf{u} \in S^{l,\rho,1,\theta}_{\beta,T}$,  the above two terms are both in $K'^{l,\rho}_{\beta,T}$. 
Since the terms in \eqref{Pbarnormbord} and \eqref{Pbartanbord} are formally identical, it is enough to prove the statement for the normal component. 
Using \eqref{inversion}, we have that
\begin{equation}
|\partial_x^i \gamma \bar{P}^\infty_n (\tilde{E}_2 \textbf{u})|_{0,\rho'} \leq  |\sup \limits_{Y \in \Gamma(\theta')}|\partial_x^i \tilde{E}_2 u_1||_{0,\rho'},
\end{equation}
where the last term is bounded by \eqref{xLinfnormal}. We also have that \begin{equation}\label{3termini}
\begin{split}
 |\partial_x^i \partial_t \gamma \bar{P}^\infty_n (\tilde{E}_2 \textbf{u})|_{0,\rho'}= |\partial_x^i  \gamma \bar{P}^\infty_n \partial_t(\tilde{E}_2 \textbf{u})|_{0,\rho'} \leq \\ \leq |\partial_x^i  \gamma \bar{P}^\infty_n \textbf{u}|_{0,\rho'}  + \varepsilon^2 | \gamma \bar{P}^\infty_n \partial_{x}^{i+2} \tilde{E}_2 \textbf{u}|_{0,\rho'}+ | \gamma \bar{P}^\infty_n \partial_x^i \partial_{YY} \tilde{E}_2 \textbf{u}|_{0,\rho'},
\end{split}
\end{equation}
with $i=0,...,l-2$. For the first term, we have \begin{equation}
\begin{split}
 |\partial_x^i  \gamma \bar{P}^\infty_n \textbf{u}|_{0,\rho'} \leq C |\sup \limits_Y | \partial_x^i   u_1||_{0,\rho'} \leq C ||\partial_x^i u_1|_{W^{1,2}(\Gamma(\theta'))}|_{0,\rho'} = \\ C ||\partial_x^i u_1|_{0,\rho'}|_{W^{1,2}(\Gamma(\theta'))} \leq C |u_1|_{l,\rho',\theta'}.
 \end{split}
\end{equation}
For the second term, we use Holder's inequality and \eqref{xLinfnormal}. For the third term, $\partial_{YY} \tilde{E}_2 \textbf{u}$ can be decomposed in two pieces: a non-boundary part, which can be treated with Holder's inequality (with exponents 1 on the exponential, $\infty$ on $\tilde{E}_2$) and with \eqref{YYLinftynormeextragamma}, and the boundary part, which can be treated with Holder inequality (this time with exponents 2 and 2) and with \eqref{YYL2normeextra}. In this last term, we have $|\varepsilon |\xi'|e^{-\varepsilon|\xi'| Y}|_{L^2_Y(\mathbb{R}^+)}= C \sqrt{\varepsilon|\xi'|}$, but the additional $|\xi'|^{1/2}$ can be balanced by the fact that we have \begin{equation}
|\xi'|^{1/2}|\gamma \partial_x^i u_1|_{0,\rho'} \leq (1 + |\xi'|) |u_1|_{l-1,\rho',1,\theta} \leq |u_1|_{l,\rho',1,\theta'}.
\end{equation} 
So the terms in \eqref{Pbarnormbord} and \eqref{Pbartanbord} are in $K^{l,\rho}_{\beta,T}$, their value is zero for $t=0$ and they can be expressed like $|\xi'| \int \limits_0^{+\infty} d y' f(\xi',y',t) k (\xi',y')$, with $|\xi'| \int \limits_0^{+ \infty} d y' |k(\xi',y')| \leq 1$: thus, all the hypothesis of proposition 3.4 in \cite{SC1998B} are satisfied.

\section{Properties of $\mathcal{S} \gamma \bar{P}^{\infty} \tilde{E}_2$}   \label{proof_prop_7punto8910}
In this appendix, we shall obtain some estimates on the operator $\mathcal{S} \gamma \bar{P}^{\infty} \tilde{E}_2$; 
combining these estimates with those obtained in appendix \ref{properties_of_PinftyE2} for $\bar{P}^\infty \tilde{E}_2$, 
we shall obtain the proof of propositions \ref{Prop2}, \ref{Prop3}, \ref{Prop4}. 

More specifically: the bounds in an $L^2_{xY}$-like setting for derivatives of order up to $l$ provide the proof of proposition \ref{Prop2}; 
the estimates of the derivatives of order $l+1$ (i.e. with an "excessive" derivative), in an $L^2_{xY}$-like setting, are used to prove proposition \ref{Prop3}; 
finally, the estimates in $L^\infty_Y L^2_x$ are used for proposition \ref{Prop4}.

For notational simplicity we shall introduce the notation
 \begin{equation}
\textbf{g}=\gamma \bar{P}^{\infty} \tilde{E}_2 \textbf{u}. 
\end{equation}
The Stokes operator, see \cite{SC1998B}, can be written explicitly as  
\begin{equation}\label{stokesoperator}
\begin{split}
\mathcal{S} \textbf{g}= \left( 
\begin{array}{ll}
- N' e^{- \varepsilon Y |\xi'|} g_n + N' (1-\bar{U}) \tilde{E}_1 V_1 \textbf{g} \\
e^{- \varepsilon|\xi'|Y} g_n + \bar{U} \tilde{E}_1 V_1 \textbf{g}
\end{array} \right)
\end{split}
\end{equation}
where  $N'$ and $V_1$ are defined as 
$$ 
N'=\frac{i \xi'}{|\xi'|}, \qquad 
V_1 \textbf{g}= g_n - N' g',
$$
$\bar{U}$ is the Ukai operator
$$
\bar{U}f= \varepsilon |\xi'| \int \limits_0^Y e^{- \varepsilon |\xi'| (Y-Y')} f(\xi',Y') d Y',
$$
and $E_1 f$ solves the heat equation on the semi-space with boundary condition $f$ and homogeneous initial datum; the explicit expression 
of $E_1$ is
$$
\tilde{E}_1 f= \int \limits_0^t ds \frac{Y}{t-s} \frac{e^{-\frac{Y^2}{4(t-s)}}}{\sqrt{4 \pi (t-s)}} \int \limits_{-\infty}^{+\infty} f(x',s) \frac{e^{-\frac{(x-x')^2}{4 \varepsilon^2(t-s)}}}{\sqrt{4 \pi (t-s) \varepsilon^2}} dx'.
$$
It is useful introduce the decomposition $\mathcal{S} \textbf{g}= \mathcal{S}^E \textbf{g}+\bar{\mathcal{S}} \textbf{g}+\mathcal{S}^* \textbf{g}$, where 
\begin{eqnarray}
\begin{split}
\mathcal{S}^E \textbf{g}= \left( 
\begin{array}{ll}
- N' e^{- \varepsilon Y |\xi'|} g_n  \\
e^{- \varepsilon|\xi'|Y} g_n 
\end{array} \right), 
\end{split}
\quad
\begin{split}
\bar{\mathcal{S}} \textbf{g}= \left( 
\begin{array}{ll}
-N' \bar{U} \tilde{E}_1 V_1 \textbf{g} \\
 \bar{U} \tilde{E}_1 V_1 \textbf{g}
\end{array} \right), 
\end{split}
\quad 
\begin{split}
\mathcal{S}^* \textbf{g}= \left( 
\begin{array}{ll}
 N'  \tilde{E}_1 V_1 \textbf{g} \\
0
\end{array} \right).
\end{split}
\nonumber
\end{eqnarray}
In what follows we shall estimate  $\mathcal{S}^E \textbf{g}$, $\bar{\mathcal{S}} \textbf{g}$ and $\mathcal{S}^* \textbf{g}$ separately. 
The general ideas to bound these terms are the following.  
The estimates of $\mathcal{S}^E \textbf{g}$ are easy to achieve, since every normal derivative is essentially transformed into $\varepsilon$ times a tangential derivative. 
Concerning $\bar{\mathcal{S}} \textbf{g}$, since all terms appear inside the Ukai operator $\bar{U}$, 
$\partial_{Y}$ operates essentially as $\varepsilon$ times a tangential derivative; 
therefore, the estimates of $\partial_Y \partial_x^j \bar{\mathcal{S}} \textbf{g}$ 
are achieved  if one is able to estimate $ \varepsilon \partial_x^{j+1} \mathcal{S}^* \textbf{g}$.  
Concerning  the second normal derivative $\partial_{YY}$, 
one has that $\partial_{YY} \partial_x^j \bar{\mathcal{S}} \textbf{g}$ has a part that behaves like 
$ \varepsilon \partial_Y \partial_x^{j+1} \mathcal{S}^* \textbf{g}$, 
and a part that behaves like $\varepsilon^2 \partial_x^{j+2} \mathcal{S}^* \textbf{g}$.
This means that the estimates of $\bar{\mathcal{S}} \textbf{g}$ follow from the estimates of $\mathcal{S}^* \textbf{g}$.

\subsection{Proof of Proposition \ref{Prop2}}
We begin with the estimates of $\mathcal{S}^E \textbf{g}$
\begin{equation}\label{SEL2}
\begin{split}
\left| \partial_x^i \partial_Y^j \mathcal{S}^E \textbf{g} \right|_{0,\rho',\theta'}^2 \leq \\ 
\varepsilon^{2j} \int \limits_{\Gamma(\theta')} d Y \int d \xi' e^{2 \rho |\xi'|} \left[ \int \limits_0^{+\infty} d Y' \varepsilon |\xi'| e^{- \varepsilon |\xi'|(Y+Y')} |\xi'|^{i+j} |\tilde{E}_2 u_1| \right]^2   \leq \\ 
\varepsilon^{2j} \int \limits_{\Gamma(\theta')} d Y \int d \xi' e^{2 \rho |\xi'|}\int \limits_0^{+\infty} d Y' \varepsilon |\xi'| e^{- \varepsilon |\xi'|(2Y+Y')} |\xi'|^{2(i+j)} |\tilde{E}_2 u_1|^2    \leq \\ 
C \varepsilon^{2j} \int d \xi' e^{2 \rho |\xi'|}\int \limits_0^{+\infty} d y'  e^{- \varepsilon |\xi'|Y'} |\xi'|^{2(i+j)} |\tilde{E}_2 u_1|^2  d y' \leq C \varepsilon^{2j} |\tilde{E}_2 u_1|_{i+j,\rho',0,\theta'}^2.
\end{split}
\end{equation}
In the second inequality we have used Cauchy-Schwartz inequality; in the third inequality we have integrated over $Y$, 
and in the last inequality we have used $|e^{- \varepsilon |\xi'|Y'}|\leq 1$. 

Now we pass to the estimates of $\mathcal{S}^* \textbf{g}$. First, we consider the case when one takes the tangential derivatives 
$\partial_{x}^j$, with $i \leq l$:  
\begin{equation}\label{SBARxL2} \begin{split}
|\partial_x^i N' \tilde{E}_1 V_1 \textbf{g}|_{0,\rho',\theta'}\leq \\
 C \int \limits_0^t ds \frac{|\gamma \bar{P}^\infty \tilde{E}_2 \textbf{u}|_{i,\rho'}}{(t-s)^{3/4}} \leq C \int \limits_0^t \frac{ds}{(t-s)^{3/4}} \int \limits_0^s  \frac{|u(\cdot, \cdot,s')|_{i,\rho',0,\theta'}}{(s-s')^{1/4}} d s' = \\  
C \int \limits_0^t |u(\cdot, \cdot,s')|_{i,\rho',0,\theta'}ds' \int \limits_{s'}^t \frac{ds}{(t-s)^{3/4} (s-s')^{1/4}} \leq  C \int \limits_0^t |u(\cdot, \cdot,s')|_{i,\rho',0,\theta'}ds' ,
\end{split}
\end{equation}
where we used the fact that 
$$
\int \limits_{s'}^t \frac{ds}{(t-s)^{3/4} (s-s')^{1/4}}= \int \limits_0^{t-s'} \frac{d s''}{(t-s' - s'')^{3/4}(s'')^{1/4}} = \int \limits_0^1 \frac{d \tau}{(1-\tau)^{3/4} \tau^{1/4}} <C.
$$

Second, one has to estimate the normal derivatives of $\mathcal{S}^* \textbf{g}$; define \begin{equation}
 f= N' V_1 \gamma \bar{P}^{\infty} \tilde{E}_2 \textbf{u}.
\end{equation}
Some lengthy but straightforward calculations lead to write 
\begin{equation}
\begin{split}
\partial_Y \tilde{E}_1 f =\int \limits_0^t ds \frac{e^{-\frac{y^2}{4\nu(t-s)}}}{\sqrt{\pi \nu (t-s)}} \int \limits_{\mathbb{R}^{d-1}}d x''  \frac{e^{- \frac{|x-x'|^2}{4 \nu (t-s)}}}{(4 \pi \nu (t-s))^{(d-1)/2}} (\nu \partial_{xx} f- \partial_t f)= \\ = \tilde{D}_1 ( \varepsilon^2 \partial_{xx}f - \partial_t f) .
\end{split}
\end{equation}

In the  functional setting of Proposition \ref{Prop2}, we are not allowed to take a time derivative: 
however, we can use 
the properties of the operator $\tilde{E}_2$ to express the time derivative as a function of the normal and tangential derivatives. 
In fact 
\begin{equation}
\partial_t f= N' V_1 \gamma \bar{P}^{\infty} \partial_t \tilde{E}_2 \textbf{u},
\end{equation}
with
\begin{equation}\label{h123}
\gamma \bar{P}^{\infty} \partial_t \tilde{E}_2 \textbf{u}= \gamma \bar{P}^\infty \textbf{u} + \gamma \bar{P}^\infty\varepsilon^2 \partial_{xx}\tilde{E}_2 \textbf{u} + \gamma \bar{P}^\infty \partial_{YY}\tilde{E}_2 \textbf{u}= \textbf{h}_1 + \textbf{h}_2 + \textbf{h}_3
\end{equation}
and 
\begin{equation}
\varepsilon^2 \partial_{xx} f= N' V_1 \varepsilon^2 \partial_{xx} \gamma \bar{P}^{\infty} \tilde{E}_2 \textbf{u},
\end{equation}
with
\begin{equation}
\varepsilon^2 \gamma \partial_{xx} \bar{P}^{\infty}  \tilde{E}_2 \textbf{u} = \textbf{h}_2.
\end{equation}

Although the contribution of $\textbf{h}_2$ is canceled, it's still useful to perform an estimate for this term, since it is used in the estimate of $\textbf{h}_3$. \\
As regards $\textbf{h}_1$, we have that, for $i \leq l-1$ \begin{equation}
|| \sup Y | \partial_x^i \textbf{u} | ||_{\rho'} \leq C || \textbf{u}||_{l,\rho',1,\theta'},
\end{equation}
with the trace which satisfies a similar bound. 
Therefore \begin{equation}
|\tilde{D}_1 N' \gamma u_1|_{i,\rho',0,\theta'} \leq C \int \limits_0^t \frac{|\textbf{u} |_{l,\rho',0,\theta'}}{(t-s)^{1/4}}ds.
\end{equation}

In order to bound $\textbf{h}_2= \gamma \bar{P}^{\infty} \varepsilon^2 \partial_{xx} \tilde{E}_2 \textbf{u}$, we let one of the tangential derivatives of $\partial_{xx} \tilde{E}_2 \textbf{u}$ hit the kernel $\tilde{E}_2$ and one hit $\textbf{u}$, so that, for $i \leq l-1$ \begin{equation}
\begin{split}
|\tilde{D}_1 V_1 \textbf{h}_2|_{i,\rho',0,\theta'} \leq C \varepsilon \int \limits_0^t \frac{ds}{(t-s)^{1/4}} \int \limits_0^s \frac{|u|_{i+1,\rho',0,\theta'}(s')}{(s-s')^{3/4}} d s' = \\ C \varepsilon \int \limits_0^t |u|_{l,\rho',0,\theta'}(s') d s'.
\end{split}
\end{equation}

We first estimate $\textbf{h}_3$, assuming that $\gamma \textbf{u}=\textbf{0}$: under this assumption, in $ \partial_{YY} \tilde{E}_2 \textbf{u}$, we can move one of the normal derivatives on $\textbf{u}$, without boundary terms. Therefore, for $i \leq l-1$ we have, in the $L^2_Y$ setting of the proposition \ref{Prop2} \begin{equation}\label{h3L2}
\begin{split}
|\tilde{D}_1 V_1 \textbf{h}_3|_{i,\rho',0,\theta'} \leq C \int \limits_0^t \frac{ds}{(t-s)^{1/4}} \int \limits_0^s \frac{\partial_Y \textbf{u}|(s')_{i,\rho',0,\theta'} }{(s-s')^{3/4}} |ds \leq C \int \limits_0^t |\textbf{u}|_{l,\rho',1,\theta'}(s) ds,
\end{split}
\end{equation}
If $\gamma \textbf{u} \neq \textbf{0}$, an additional term appears in the expression of $\bar{P}^\infty \partial_{YY} \tilde{E}_2 \textbf{u}$, which roughly behaves like $\tilde{E}_1 \gamma u_1$; the loss of half derivative due to the trace operator is not a problem for the tangential derivatives of order up to $l-1$, so inequality \eqref{h3L2}  still holds.
This concludes the estimates of the normal derivative of $\bar{\mathcal{S}}$.
 \\
As we already said at the beginning of this section, the estimates of $\bar{\mathcal{S}} \textbf{g}$ can be derived from the ones of $\mathcal{S}^* \textbf{g}$: therefore, the proof of the proposition is complete.

\subsection{Proof of Proposition \ref{Prop3}}
The estimate \eqref{SEL2} obtained for $\mathcal{S}^E \textbf{g}$ still holds for derivatives of order $l+1$; the "excessive" derivative can be treated either using the Cauchy estimate, or giving a tangential derivative to the kernel of $\tilde{E}_2$.

The same applies for $\partial_x^{l+1} \mathcal{S}^* \textbf{g}$: the only difference is that we change the order of integration in time. This is done, when we use the Cauchy estimate, in order to eliminate the mild singularity in time deriving from the $L^2_Y$ norm of the operato $\tilde{E}_1$, while, when we let one tangential derivative hit the kernel of $\tilde{E}_2$, changing the order of integration allows to reduce the singularity in time. We have \begin{equation} \label{SBARxL2extra}
|N' \tilde{E}_1 V_1 \textbf{g}|_{l+1,\rho',0,\theta'} \leq \frac{C}{\varepsilon} \int \limits_0^t |u|_{l,\rho',0,\theta'}ds' \int \limits_{s'}^t \frac{ds}{(t-s)^{3/4} (s-s')^{3/4}},
\end{equation}
with \begin{equation} \begin{split}
\int \limits_{s'}^t \frac{ds}{(t-s)^{3/4} (s-s')^{3/4}}= \int \limits_0^{t-s'} \frac{d s''}{(t-s - s'')^{3/4}(s-s')^{3/4}} \\= \frac{1}{(t-s')^{1/2}} \int \limits_0^1 \frac{d \tau}{(1-\tau)^{3/4} \tau^{3/4}},
\end{split}
\end{equation}
so \begin{equation}\label{SBARxL2extrafinale}
|N' \tilde{E}_1 V_1 \textbf{g}|_{l+1,\rho',0,\theta'} \leq \frac{C}{\varepsilon} \int \limits_0^t \frac{|u|_{l,\rho',0,\theta'}(s')}{(t-s')^{1/2}}ds' .
\end{equation}
\\
Now, we begin to estimate $\partial_Y \mathcal{S}* \textbf{g}$.
First, we estimate the term $\textbf{h}_1$, introduced in \eqref{h123}, under the assumption that $\gamma \textbf{u}=0$: in this case,  \begin{equation}
|\textbf{h}_1| \leq C \sqrt{\varepsilon |\xi'|} |\textbf{u}|_{L^2_Y}.
\end{equation}
Using \begin{equation}\label{radicesukernelE1}
(\varepsilon |\xi'|)^{1/2} e^{- \varepsilon^2 (t-s) |\xi'|^2} \leq \frac{C}{(t-s)^{1/4}}
\end{equation}
we have that \begin{equation}\label{SBARYpt1L2}
|\tilde{D}_1 V_1 \textbf{h}_1|_{l,\rho',0,\theta'} \leq C \int \limits_0^t \frac{|\textbf{u} |_{l,\rho',0,\theta'}}{(t-s)^{1/2}}ds.
\end{equation}
If $\gamma \textbf{u} \neq 0$, then $V_1 \textbf{h}_1$ contains $N' \gamma u_1$ as additional term, with \begin{equation}\label{interp}
|\gamma u_1| \leq |u_1|_{L^2_Y}^{1/2} |\partial_Y u_1|_{L^2_Y}^{1/2}.
\end{equation}
We use the interpolation inequality \eqref{interp}, together with the Cauchy estimates, to obtain \begin{equation}
|\tilde{D}_1 N' \gamma u_1|_{l,\rho',0,\theta'} \leq C \int \limits_0^t \frac{|\textbf{u} |^{1/2}_{l,\rho',0,\theta'}|\textbf{u} |^{1/2}_{l,\rho(s),1,\theta'}}{(t-s)^{1/4}(\rho(s)-\rho')^{1/2}}ds \leq C \int \limits_0^t \frac{|\textbf{u} |_{l,\rho',0,\theta'}|}{(t-s)^{1/2}} + \frac{|\textbf{u} |_{l,\rho(s),1,\theta'}}{(\rho(s)-\rho')}ds.
\end{equation}
This concludes the estimates of $\textbf{h}_1$.
\\ The bounds of $\textbf{h}_2= \gamma \bar{P}^{\infty} \varepsilon^2 \partial_{xx} \tilde{E}_2 \textbf{u}$ are obtained moving one tangential derivative on the kernel of $D_1$ and one tangential derivative on the kernel of $\tilde{E}_2$, so \begin{equation}
|\tilde{D}_1 V_1 \textbf{h}_2|_{l,\rho',0,\theta'} \leq C  \int \limits_0^t \frac{ds}{(t-s)^{3/4}} \int \limits_0^s \frac{|u|_{l,\rho',0,\theta'}(s')}{(s-s')^{3/4}} d s' = C  \int \limits_0^t \frac{|u|_{l,\rho',0,\theta'}(s')}{(t-s')^{1/2}} d s'.
\end{equation}
 \\
We begin the estimates of $\textbf{h}_3$, starting from the case $\gamma \textbf{u}=\textbf{0}$: in this case, $\gamma \partial_YY \tilde{E}_2 \textbf{u} =\textbf{0}$, so only the integral part of $\gamma \bar{P}^\infty$ is nonzero, when we apply this operator to $\partial_YY \tilde{E}_2 \textbf{u}$. We can integrate by parts two times, moving the normal derivatives on the kernel of $\bar{P}^\infty$: the result is a term which behaves like $\varepsilon^2 \partial_{xx} \bar{P}^\infty \tilde{E}_2 \textbf{u}$ (and therefore can be bounded exactly like $\textbf{h}_2$) plus an additional term, given by \begin{equation}
C i \xi' \varepsilon \int \limits_0^s d s' e^{- \varepsilon^2 (s-s') \xi^{\prime 2}} \int \limits_0^{+ \infty} \frac{Y'}{(s-s')^{3/2}} e^{-\frac{Y'^2}{4(s-s')}} (u_2,u_1)(\xi',Y,s) d Y'.
\end{equation}
This term is inside $D_1$: we can move $i \xi'$ on the kernel of $D_1$, so that the singularities in time are better distributed. Therefore, we obtain \begin{equation}
|\tilde{D}_1 V_1 \textbf{h}_3|_{l,\rho',0,\theta'} \leq  C  \int \limits_0^t \frac{|u|_{l,\rho',0,\theta'}(s')}{(t-s')^{1/2}} d s'.
\end{equation}
If $\gamma \textbf{u} \neq \textbf{0}$, an additional term appears in the expression of $\bar{P}^\infty \partial_{YY} \tilde{E}_2 \textbf{u}$, which roughly behaves like $\tilde{E}_1 \gamma u_1$; as for $\textbf{h}_1$, in this case we use the interpolation inequality \eqref{interp}, the Cauchy estimate and Young's inequality, obtaining \begin{equation}
|\tilde{D}_1 V_1 \textbf{h}_3|_{l,\rho',0,\theta'} \leq  C  \int \limits_0^t \frac{|u|_{l,\rho',0,\theta'}}{(t-s')^{1/2}} + \frac{|u|_{l,\rho(s'),1,\theta'}(s')}{\rho(s')- \rho'} ds'.
\end{equation}
\\
Finally, when we take another derivative with respect to $Y$, we obtain, for a generic $f$ with $f(t=0)=0$, that \begin{equation}
\begin{split}
\partial_{YY} \tilde{E}_1 f= \\ \int \limits_0^t \frac{Ye^{-\frac{Y^2}{4(t-s)}}}{4 \sqrt{\pi} (t-s)^{3/2}} ds \int \limits_{- \infty}^{+ \infty} \left((\partial_t f)(x',s)- \varepsilon^2 (\partial_{xx}f)(x',s) \right) \frac{e^{-\frac{(x-x')^2}{4 \varepsilon^2 (t-s)}}}{\sqrt{4 \pi \varepsilon^2 (t-s)}} dx'  \\ = \tilde{E}_1 (\partial_t f - \varepsilon^2 \partial_{xx}f ).
\end{split}
\end{equation} 
Therefore, we have to evaluate $\tilde{E}_1 \gamma V_1 \textbf{h}_1$ and $ \tilde{E}_1 \gamma V_1 \textbf{h}_3$: this time, since two normal derivatives have been taken, only up to $l-1$ tangential derivatives can be taken. Therefore, this time, the loss of half derivative due to the trace operator is not a problem, and the case $\gamma \textbf{u} \neq \textbf{0}$ does not need to be treated separately. For $\textbf{h}_1$, we have  \begin{equation}
|\tilde{E}_1 V_1 \textbf{h}_1|_{i,\rho',0,\theta'} \leq  C  \int \limits_0^t \frac{|\textbf{u}|_{i+i,\rho',0,\theta'}}{(t-s)^{3/4}} d s.
\end{equation} 
For $\tilde{E}_1 \gamma V_1 \textbf{h}_3$, a similar estimate holds: therefore, the estimates of $\mathcal{S}^*$ are complete.
\\ The estimates of $\bar{\mathcal{S}} \textbf{g}$ can be derived from the ones of $\mathcal{S}^* \textbf{g}$.
\subsection{Proof of proposition \ref{Prop4}}
In the $L^\infty_YL^2_x$-like setting of proposition \ref{Prop4}, the estimates differ from the ones of \ref{Prop3} an additional $(t-s)^{-1/4}$, which implies a stronger (but still integrable) singularity in time. Indeed, for the estimates of $\mathcal{S}^* \textbf{g}$ (and therefore of $\bar{\mathcal{S}} \textbf{g}$, we give the $L^\infty_Y$ norm to the kernel of $\textbf{E}_1$, which is "$(t-s)^{-1/4}$ times worse" of its $L^2_Y$ norm. For the estimates of $\mathcal{S}^E \textbf{g}$, instead, we give $e^{- \varepsilon |\xi'| Y}$ the $L^\infty_Y$ norm, which is $1$, while its $L^2_Y$ norm is $C (\varepsilon |\xi'|)^{-1/2}$: moving the additional $(\varepsilon |\xi'|)^{1/2}$ on the kernel of $\tilde{E}_2$, we can use the bound given by equation \eqref{radicesukernelE1}, and an additional $(t-s)^{-1/4}$ appears.

\section{Estimate of the linear terms in \eqref{def::F}}\label{sub::est_linear_F}

In $ [ \textbf{u}^{NS0} + \varepsilon (\boldsymbol{w}+ \boldsymbol{w}+ \boldsymbol{\sigma} )] \cdot \nabla \mathcal{N}^* \textbf{e}^* $ the most problematic term is $u^S \partial_x \mathcal{N}^* \textbf{e}^*$: using propositions \ref{Prop2} and \ref{Prop3}, we have that \begin{equation}
|u^S \partial_x \mathcal{N}^* \textbf{e}^*|_{l,\rho',0,\theta'} \leq C \int \limits_0^t \frac{|\textbf{e}^*|_{l,\rho(s),0,\theta}}{\rho(s)-\rho'} ds,
\end{equation}
\begin{equation}
|\partial_Y u^S \partial_x \mathcal{N}^* \textbf{e}^*|_{l-1,\rho',0,\theta'} \leq \frac{C}{\sqrt{t}} \int \limits_0^t |\textbf{e}^*|_{l,\rho',0,\theta}ds,
\end{equation}
\begin{equation}
| u^S \partial_Y \partial_x \mathcal{N}^* \textbf{e}^*|_{l-1,\rho',0,\theta'} \leq C \int \limits_0^t \frac{|\textbf{e}^*|_{l,\rho',1,\theta}}{(t-s)^{1/2}}ds.
\end{equation}
In $\mathcal{N}^* \textbf{e}^* \cdot \nabla [ \textbf{u}^{NS0}  + \varepsilon (\boldsymbol{\sigma} + \boldsymbol{w} + \boldsymbol {\sigma})] $, the most problematic term is $\frac{1}{\varepsilon}(\mathcal{N}^* \textbf{e}^*)_2 \partial_Y \tilde{u}^P$ (the problem with $\tilde{u}^R$ is to prove that the term is $O(1)$). Since $ \mathcal{N}^* \textbf{e}^*$ is divergence free and zero at the boundary, we have \begin{equation}
\left( \mathcal{N}^* \textbf{e}^* \right)_2 = - \varepsilon \int \limits_0^Y \partial_x \left( \mathcal{N}^* \textbf{e}^* \right)_1 d Y',
\end{equation}
so \begin{equation}
\begin{split}
\left| \frac{1}{\varepsilon}(\mathcal{N}^* \textbf{e}^*)_2 \partial_Y \tilde{u}^P \right|_{l,\rho',0,\theta'} \leq C  \left| \frac{1}{Y}\int \limits_0^Y \partial_x \left( \mathcal{N}^* \textbf{e}^* \right)_1 d Y' \right|_{l,\rho',0,\theta'} \sup \limits_Y \left| Y \partial_Y \tilde{u}^P \right|_{l,\rho'} \leq \\ C \left| \partial_x \left( \mathcal{N}^* \textbf{e}^* \right)_1 \right|_{l,\rho',0,\theta'} \leq C \int \limits_0^t \frac{|\textbf{e}^*|_{l,\rho(s),0,\theta'}}{\rho(s)-\rho'} d s.
\end{split}
\end{equation}
In the first inequality we used algebra property \ref{algebrax}, in the second inequality we used the boundedness of $\sup \limits_Y \left| Y \partial_Y \tilde{u}^P \right|_{l,\rho'}$ (both for the regular and for the singular part) and the fact that the average operator is bounded from $L^2$ to $L^2$ (is bounded by the maximal function, so its trivial), and in the last inequality we used proposition \ref{Prop4}. With a similar argument we have that \begin{equation}
\left| \frac{1}{\varepsilon}(\mathcal{N}^* \textbf{e}^*)_2 \partial_{YY} \tilde{u}^P \right|_{l-1,\rho',0,\theta'} \leq C \left( 1 + \frac{1}{t^{1/2}} \right) \int \limits_0^t |\textbf{e}^*|_{l,\rho',0,\theta'} ds
\end{equation}
and \begin{equation}\label{ultimaquasicontr}
\begin{split}
\left| \frac{1}{\varepsilon}\partial_Y (\mathcal{N}^* \textbf{e}^*)_2 \partial_{Y} \tilde{u}^P \right|_{l-1,\rho',0,\theta'}= \left| \partial_x(\mathcal{N}^* \textbf{e}^*)_1 \partial_{Y} \tilde{u}^P \right|_{l-1,\rho',0,\theta'} \leq \\  C \left( 1 + \frac{1}{t^{1/2}} \right) \int \limits_0^t |\textbf{e}^*|_{l,\rho',0,\theta'} ds
\end{split}.
\end{equation}

\bibliographystyle{siamplain}
\bibliography{references} 
\end{document}

%% file: ex_shared.tex
\usepackage{lipsum}

\usepackage{graphicx,color}

\usepackage{algorithmic}

\title{Navier-Stokes equations in the half space with non compatible data}

\author{Andrea Argenziano \thanks{Dipartimento di Matematica e Informatica, Universit\`a di Palermo, Via Archirafi 34, 90123 Palermo, Italy
  (\email{andrea.argenziano@unipa.it}).}
\and Marco Cannone \thanks{Laboratoire d'Analyse et de Mathématiques Appliquées
Université Gustave Eiffel
5 boulevard Descartes
Bâtiment Copernic
77420 Champs-sur-Marne
  (\email{marco.cannone@univ-eiffel.fr}).}
\and Marco Sammartino\thanks{Dipartimento di Ingegneria, Universit\`a di Palermo, Viale delle Scienze Ed.8, 90100 Palermo, Italy
  (\email{marco.sammartino@unipa.it}).}}

\usepackage{amsopn}

\usepackage{amsfonts,amsmath,amsthm}

\theoremstyle{plain}% default
\newtheorem{theorem}{Theorem}[section]
\newtheorem{proposition}{Proposition}[section]
\newtheorem{lemma}{Lemma}[section]

\newtheorem*{keywords}{Keywords}
\newtheorem*{AMS}{AMS}

\theoremstyle{definition}
\newtheorem{definition}{Definition}[section]
\newtheorem{remark}{Remark}[section]

\numberwithin{equation}{section}

\usepackage{hyperref}
\newcommand*{\email}[1]{\href{mailto:#1}{\nolinkurl{#1}}}
\hypersetup{
  pdftitle={Navier Stokes equations in the half space with non compatible data},
  pdfauthor={M. Cannone, A. Argenziano, and M. Sammartino}
}

\usepackage{cleveref}

\def\color#1{\relax}